%% file: compactart6.tex
\documentclass[11pt]{smfart}

\input{paquets.tex}
\input{numerotation}

\input{lemme-proposition-theoreme.tex}

\numberwithin{equation}{subsubsection} 
\renewcommand{\theequation}{\thesubsubsection.\Alph{equation}}

\newcommand{\agb}{\overline{\mathcal{A}}_{g}}
\newcommand{\agKb}{\overline{\mathcal{A}}_{g,\, 0}}

\newcommand{\iso}{\buildrel{\sim}\over\rightarrow}
\newcommand{\isolong}{\buildrel{\sim}\over\longrightarrow}

\newcommand{\Z}{\mathbb Z}
\newcommand{\Q}{\mathbb Q}
\newcommand{\R}{\mathbb R}

\newcommand{\F}{\mathbb F}

\newcommand{\Gm}{{\mathbb G}_m}
\newcommand{\Mot}{\mathsf{Mot}}
\newcommand{\Mum}{\mathsf{Mum}}
\newcommand{\RHom}{\underline{\mathrm{RHom}}}
\newcommand{\D}{\mathcal{D}^{b}}
\newcommand{\K}{\mathcal{C}^{b}}
\newcommand{\Y}{\mathcal{Y}}
\newcommand{\gsp}{\mathrm{GSp}_{2g}}
\newcommand{\ag}{\mathcal{A}_{g}}
\newcommand{\agK}{\mathcal{A}_{g , \, 0}}
\newcommand{\Ag}{{\mathrm{A}}_{g}}

\newcommand{\agr}{\mathcal{A}_{g-r}}
\newcommand{\agrK}{\mathcal{A}_{g-r , \, 0}}
\newcommand{\Agr}{{\mathrm{A}}_{g-r}}

\newcommand{\Xsi}{\mathcal{M}_{\sigma}}
\newcommand{\XsiK}{\mathcal{M}_{\sigma,\,0}}
\newcommand{\XsiKb}{\overline{{\mathcal{M}}}_{\sigma,\,0}}
\newcommand{\Xsiw}{\mathcal{M}^{w}_{\sigma,\, 0}}
\newcommand{\Xsib}{\overline{\mathcal{M}}_{\sigma}}
\newcommand{\Xsiwb}{\overline{\mathcal{M}}^{w}_{\sigma, \ 0}}

\newcommand{\Hom}{\underline{\mathrm{Hom}}}
\newcommand{\Ext}{\underline{\mathrm{Ext}}^1}

\newcommand{\Spec}{\mathrm{Spec}}
\def\sga#1#2#3{[{\bf $\mathbf{SGA\,{#1}}$}~{\sc #2}~#3]} 	
\def\ega#1#2{[{\bf \'EGA}~{\sc #1}~#2]} 				
\def\egalong#1#2#3{[{\bf \'EGA}~$\textsc{#1}_{\textrm{#2}}$~#3]} 	
\def\egazéro#1#2{[{\bf \'EGA}~$0_{\textsc{#1}}$~#2]}		

\title{Compactification de variétés de Siegel aux places de mauvaise réduction}
\author{Benoît Stroh}
\date{17 novembre 2008}
\email{benoit.stroh@gmail.com}
\address{Institut Élie Cartan, Université Henri Poincaré,  B.P.~239, F-54506 Vandoeuvre-les-Nancy Cedex, France}

\begin{document}

\begin{abstract}
Nous construisons des compactifications toroïdales arithmétiques du champ de modules des variétés abéliennes principalement polarisées munies d'une structure de niveau parahorique. Pour ce faire, nous étendons la méthode de Faltings et Chai~\cite{Deg@FaltingsChai} à un cas de mauvaise réduction. Le voisinage du bord des compactifications obtenues n'est pas lisse, mais a pour singularités celles des champs de modules de variétés abéliennes avec structure parahorique de genre plus petit. Nous sommes amenés à reprendre la construction des compactifications sans niveau de Faltings et Chai, en modifiant l'étape d'approximation pour préserver le groupe de $p$-torsion des variétés abéliennes. Nous donnons comme application une nouvelle preuve de l'existence du sous-groupe canonique pour des familles de variétés abéliennes.
\end{abstract}

\begin{altabstract} \textbf{(Compactification of Siegel modular varieties with bad reduction)}\\ 
We construct arithmetic toroidal compactifications of the moduli stack of principally polarized abelian varieties with parahoric level structure. To this end, we extend the methods of Faltings and Chai~\cite{Deg@FaltingsChai} to a case of bad reduction. Our compactifications are not smooth near the boundary; the singularities are those of the moduli stacks of abelian varieties with parahoric level structure of lower genus. We modify Faltings and Chai's construction of compactifications without level structure. The key point is that our approximation preserves the $p$-torsion subgroup of the abelian varieties. As an application, we give a new proof of the existence of the canonical subgroup for some families of abelian varieties.
\end{altabstract}

\maketitle

\newpage

{ \small \tableofcontents}

Dans cet l'article, $g$ désigne un entier naturel et l'on note $\ag$ le champ de modules des schémas abéliens principalement polarisés de genre $g$. Ce champ est algébrique et lisse sur $\Spec(\Z)$, mais pas propre. En effet, le critère valuatif de propreté n'est pas vérifié car si ${R}$ est un anneau de valuation discrète de corps des fractions~$K$, il existe des variétés abéliennes principalement polarisées sur $\Spec({K})$ qui ne se prolongent pas en des schémas abéliens sur $\Spec(R)$. Il est donc naturel de construire des compactifications de $\ag$. Ce problème, initié par Mumford, est riche d'une longue histoire~\cite{ReviewFC@Raynaud}. Ash, Mumford, Rapoport et Tai~\cite{Compact@AMRT} d'une part, Faltings et Chai~\cite{Deg@FaltingsChai} d'autre part, ont construit des \textit{compactifications toroïdales} de $\ag$. Ces dernières sont des champs algébriques dépendant du choix d'un éventail $\Sigma$. Notons 
$${\overline{\mathcal{A}}_g^{\Sigma}}$$
la compactification toroïdale associée à un tel éventail $\Sigma$. Sous une hypothèse de nature combinatoire sur $\Sigma$, le champ $${\overline{\mathcal{A}}_g^{\Sigma}}$$ est lisse sur $\Spec(\Z)$ et le complémentaire de $\ag$ forme un diviseur à croisements normaux. La stratification induite par ce diviseur est paramétrée par $\Sigma$ et l'on peut associer à chaque strate un rang $r\leq g$. Le complété formel de la compactification le long d'une strate de rang $r$ est isomorphe au complété formel d'un certain fibré en plongements toriques $\overline{\mathcal{M}}^r \rightarrow \mathcal{B}^r$
le long d'une strate, où $\mathcal{B}^r \rightarrow \agr$ est un schéma abélien.

Si $g\geq 2$, on ne dispose malheureusement pas d'interprétation modulaire des compactifications toroïdales de $\ag$. On peut toutefois affirmer qu'elles paramètrent 
``tous les schémas semi-abéliens sur les traits qui sont génériquement des schémas abéliens principalement polarisés''. La propreté de ces compactifications résulte alors du critère valuatif et du théorème de réduction semi-stable pour les variétés abéliennes~\sga{7}{ix}{th.~3.6}.

\subsection*{\'Enoncé du théorème principal}

Dans la suite de l'article, on fixe un nombre premier $p$, un entier $s\leq g$ et un ensemble ordonné non vide $$\mathcal{D}=\{d_1<d_2<\cdots<d_s\}\subset\{1,\cdots,g\}\: .$$
Définissons le champ $\agK$ comme l'espace de modules relatif sur $\ag$ qui paramètre les drapeaux de sous-groupes finis, plats et totalement isotropes
$$H_1 \: \subset \: H_2 \: \subset \: \cdots \: \subset \: H_s $$
du groupe de $p$-torsion de la variété abélienne universelle, tels que 
$H_i$ soit de rang $d_i$ pour tout $i\leq s$. On dit que $\agK$ est le champ des schémas abéliens principalement polarisés de genre $g$ munis d'une structure de niveau parahorique de type $\mathcal{D}$ en $p$. Ce champ est algébrique mais n'est pas lisse sur $\Spec(\Z_p)$~;~on dit qu'il a \textit{mauvaise réduction} en $p$. Il n'est pas propre sur $\Spec(\Z)$, et le but de cet article est d'en construire des compactifications toroïdales. Bien sûr, on ne peut pas obtenir de compactifications lisses~;~on voudrait toutefois contrôler leurs singularités. Le résultat suivant sera précisé dans le théorème~\ref{thPrinc1}.

\begin{thprinc} Il existe un champ algébrique canonique 
${\overline{\mathcal{A}}_{g,\: 0}^{\Sigma}}$
qui est propre sur $\Spec(\Z)$, qui contient $\agK$ comme ouvert dense et tel qu'il existe une flèche propre $${\overline{\mathcal{A}}_{g,\: 0}^{\Sigma}}\longrightarrow {\overline{\mathcal{A}}_g^{\Sigma}}$$ prolongeant le morphisme d'oubli du niveau. Le champ ${\overline{\mathcal{A}}_{g,\: 0}^{\Sigma}}$ est muni d'une stratification paramétrée par~$\Sigma$ et l'on peut décrire son complété formel le long d'une strate par certains fibrés en plongements toriques. En particulier, ses singularités sont produits de singularités toriques et des singularités des champs $\agrK$ pour $r\geq 0$.
\end{thprinc}

La présence de singularités toriques est due au caractère trop restrictif du choix combinatoire $\Sigma$. Pour remédier à ce problème, nous remplaçons le choix de $\Sigma$ par celui d'un éventail~$\mathfrak{S}$ recouvrant un complexe conique adéquat. Dans le théorème~\ref{thPrinc}, nous construisons une compactification toroïdale \og améliorée \fg dépendant de $\mathfrak{S}$. Ses singularités sont celles des champs $\agrK$ pour $r\geq 0$. En particulier, la compactification améliorée est lisse sur $\Spec(\Z[1/p])$ et le complémentaire de $\agK$ est un diviseur à croisements normaux.

La méthode de Faltings et Chai utilise de manière cruciale la lissité de $\ag$, et n'est donc pas directement transposable aux cas de mauvaise réduction. Nous contournons ce problème en reprenant la construction des compactifications de $\ag$ et en la modifiant. En particulier, notre approximation sur la variété sans niveau \textit{préserve le sous-groupe de $p$-torsion}.
On trouvera donc dans cette thèse la construction simultanée des compactifications de $\ag$ et de $\agK$. Nous donnons une présentation différente des résultats de~\cite{Deg@FaltingsChai} en introduisant les notions de {1-motif de Mumford} et de {catégorie fibrée à la frontière}.

La technique utilisée ici est transposable à de multiples cas~P.E.L, comme celui des variétés de Hilbert-Siegel ou celui des variétés de Shimura unitaires, en utilisant les résultats de~\cite{FlatSympl@Gortz} et de~\cite{FlatUnitary@Gortz}. Toutefois, la construction de compactifications toroïdales en une place de mauvaise réduction qui divise le degré de la polarisation ou le discriminant de l'anneau d'endomorphismes reste un problème ouvert~;~les méthodes de~\cite{Deg@FaltingsChai} ne semblent pas s'y adapter.

Nous donnons comme application une nouvelle preuve de l'existence du sous-groupe canonique pour des familles de variétés abéliennes, due à Abbès et Mokrane~\cite{Canonique@AbbesMokrane} et à Andreatta et Gasbarri~\cite{Canonique@Andreatta}. Nous utilisons une idée due à Vincent Pilloni. Supposons que $\mathcal{D}=\{g\}$~;~on parle alors de \textit{structure de niveau de Siegel}. Soit $\agKb$ (resp.~$\agb$) une compactification toroïdale de $\agK$ (resp. de~$\ag$). Nous définissons le lieu ordinaire $\agb^{\mathrm{ord}}$ de $\agb\times \Spec(\F_p)$ et le lieu ordinaire-multiplicatif $\agKb^{\mathrm{ord-mult}}$ de $\agKb\times \Spec(\F_p)$. Nous montrons que la section de $\ag^{\mathrm{ord}}$ dans $\agK^{\mathrm{ord-mult}}$ donnée par le noyau du morphisme de Frobenius s'étend en une section
$$\agb^{\mathrm{ord}}\:\longrightarrow \:\agKb^{\mathrm{ord-mult}}\: .$$
Nous en déduisons le théorème suivant grâce à la \og surconvergence du tube \fg, démontrée dans~\cite[th.~1.3.5]{Rigide@Berthelot}. 

\begin{thrien} La section de $\agb^{\mathrm{ord}}$ dans $\agKb^{\mathrm{ord-mult}}$ induit un isomorphisme entre des voisinages stricts des tubes de chacun de ces deux schémas dans les espaces rigides associés à $\agb$ et à $\agKb$. 
\end{thrien}

Ce théorème montre l'existence du sous-groupe canonique sur un voisinage strict du tube du lieu ordinaire. Notre méthode ne permet pas de préciser quantitativement ce voisinage, contrairement à~\cite{Canonique@AbbesMokrane} et à~\cite{Canonique@Andreatta}.

\subsection*{Plan de l'article}  Dans la première partie, on définit la notion de groupe de $p$-torsion d'un $1$-motif et celle de structure de niveau parahorique sur un $1$-motif. On traite de la construction de Mumford et l'on introduit la catégorie des  1-motifs de Mumford, qui est une version plus synthétique de la catégorie $\mathrm{DD}_{pol}$ de Faltings et Chai. On décrit également les espaces de modules de $1$-motifs avec et sans niveau.

Dans la seconde partie, nous construisons les cartes locales formelles des compactifications. Nous compactifions partiellement de manière torique les espaces de modules de $1$-motifs, en utilisant leur  description comme torseurs sous des tores~;~c'est ici qu'intervient le choix combinatoire $\Sigma$. Nous complétons ces compactifications partielles le long de strates, pour pouvoir effectuer une construction de Mumford. On remarque que le groupe de $p$-torsion du quotient à la Mumford est une donnée {algébrique}, qui existe sur l'espace de modules des $1$-motifs avant complétion~;~cela s'avérera essentiel dans la suite. 
Nous approximons ensuite le résultat de la construction de Mumford. C'est l'étape la plus délicate de la construction sans niveau, qui requiert la {lissité} des espaces de modules. En effet, il faut s'assurer que l'approximation n'entraîne pas de perte d'information (proposition~\ref{propRecomplAutoEt}). Ce contrôle s'effectue grâce à l'isomorphisme de Kodaira-Spencer (proposition~\ref{propMumUniv}), qui n'existe que sur des bases lisses. Nous  construisons une approximation préservant le groupe de $p$-torsion (proposition~\ref{propGxih})~;~cela a un sens puisque ce groupe est une donnée algébrique. 
Nous définissons enfin les cartes locales algébriques. Celles relatives à $\agK$  n'auraient aucun intérêt si le groupe de $p$-torsion n'avait pas été préservé lors de l'approximation.

Dans la troisième partie, nous construisons les compactifications toroïdales. Nous recollons les cartes locales algébriques grâce à la construction de relations d'équivalence étales. Nous utilisons à nouveau la lissité dans la construction sans niveau \textit{via} l'isomorphisme de Kodaira-Spencer (proposition~\ref{propUagEt}).  Comme $\agKb$ est propre sur $\agb$, il suffit de vérifier la propreté de la compactification sans niveau (prop.~\ref{propPROPRE}). Cela repose sur le {contrôle de l'approximation} permis par l'isomorphisme de Kodaira-Spencer.
Nous construisons ensuite une compactification \og améliorée \fg, qui n'a pas de singularités d'origine torique. Pour cela, nous remplaçons le choix de $\Sigma$ par celui d'une décomposition polyédrale $\mathfrak{S}$ d'un certain complexe conique. De plus, la combinatoire de la stratification de $\agKb$ se lit bien plus aisément sur $\mathfrak{S}$ que sur $\Sigma$.

Dans la quatrième partie, nous redémontrons l'existence du sous-groupe canonique en suivant une idée de Vincent Pilloni. Nous y utilisons un théorème de \og surconvergence du tube \fg de~\cite{Rigide@Berthelot}.

\medskip

Je remercie chaleureusement Alain Genestier, qui m'a initié aux compactifications et à la mauvaise réduction des variétés de Shimura, de m'avoir consacré tant de temps lors de longues discussions mathématiques. Je suis également redevable à Luc Illusie de la démonstration du lemme~\ref{lemApprox}, et à Vincent Pilloni de l'idée utilisée dans la démonstration du théorème~\ref{theo}.

\section{Variétés abéliennes à réduction semi-stable et $1$-motifs}

\subsection{Schémas semi-abéliens} Soit $S$ un champ algébrique.

\begin{definition2} Un schéma semi-abélien $G\rightarrow S$ est un schéma en groupes lisse, séparé, à fibres géométriquement connexes,  tel que pour tout $s\in S$ la fibre $G_s$ soit extension d'une variété abélienne $A_s$ par un tore $T_s$.
\end{definition2}

\begin{theoreme2}[\cite{Deg@FaltingsChai}, proposition~I.2.7]\label{thRay}  Soient $S$ un champ algébrique n{\oe}thérien normal, $U$ un ouvert dense de $S$, et $G$, $G'$ deux schémas semi-abéliens sur $S$. Le morphisme de restriction suivant induit un isomorphisme
$$\mathrm{Hom}_S(G,G') \isolong \mathrm{Hom}_U(G_U,G'_U).$$
\end{theoreme2}

Soient $S$ un champ algébrique n{\oe}thérien et $G\rightarrow S$ un schéma semi-abélien. Pour tout $s\in S$, le tore $T_s$ définit un faisceau étale en $\Z$-modules libres de rang fini sur $s$ en posant
$${\underline{X}}_s = {\Hom}_s(T_s , \Gm).$$
Le théorème suivant montre que ces constructions se recollent.

\begin{theoreme2}[\cite{Deg@FaltingsChai}, théorème I.2.10] \label{thFaisCons}  Il existe un unique faisceau étale constructible $\underline{X}$ sur $S$ tel que pour tout  point $s$ de $S$, on ait $\underline{X} \times_S s = \underline{X}_s$.
\end{theoreme2}

On dit que $G$ est de \textit{rang torique constant} si le faisceau étale~$\underline{X}$ est localement constant sur $S$.

\begin{proposition2}[\cite{Deg@FaltingsChai}, corollaire~I.2.11]  Soient $S$ un schéma n{\oe}thérien et $G\rightarrow S$ un schéma semi-abélien de rang torique constant. Il existe un tore $T\rightarrow S$, un schéma abélien $A\rightarrow S$ et une manière d'écrire $G$ comme extension sur $S$
$$0\longrightarrow T \longrightarrow G \longrightarrow A\longrightarrow 0.$$
\end{proposition2}

Dans la suite de cet article les schémas semi-abéliens de rang torique constant seront notés avec un \textit{tilde}, comme par exemple $\tilde{G}$. Les classes d'isomorphismes d'extensions $$0\longrightarrow T \longrightarrow \tilde{G} \longrightarrow A\longrightarrow 0$$ sont paramétrées par $\Ext_{\Z}(A,T)$, qui est égal à $\Hom_{\Z}(\underline{X},A^t)$ par la formule de Barsotti-Weil.

\subsection{Les $1$-motifs}

On se donne un champ algébrique $S$ et l'on note $\K_S$ la catégorie des complexes bornés de faisceaux \textit{fppf} en groupes abéliens sur $S$ et $\D_S$ la catégorie dérivée associée. On dit qu'un faisceau étale sur $S$ est un \textit{réseau} s'il est localement constant à valeurs dans un $\Z$-module libre de rang fini.

\begin{definition2}[\cite{HodgeIII@Deligne}, définition 10.1.1]  Un {$1$-motif} $M\rightarrow S$ est la donnée d'un schéma semi-abélien de rang torique constant $\tilde{G} \rightarrow S$,
d'un réseau $Y \rightarrow S$ et d'un complexe de $\K_S$ concentré en degrés $-1$ et $0$ $$M=[Y\rightarrow \tilde{G}].$$
\end{definition2}

\begin{remarque2} Comme tous les quasi-automorphismes de ce complexe  sont des automorphismes (\cite{1Mot@Raynaud},~10.2.14), on peut indifféremment considérer $M$ comme un objet de $\K_S$ ou de $\D_S$.
\end{remarque2}

On note $\Mot_S$ la catégorie des $1$-motifs sur $S$. Le \textit{genre} d'un $1$-motif $M$ est l'entier
$$\mathrm{rg}_{\Z}(Y)+\mathrm{dim}_S \tilde{G}\: .$$
Soit $0\rightarrow T \rightarrow \tilde{G} \rightarrow A \rightarrow 0$ l'extension qui définit $\tilde{G}$ sur $S$. Le $1$-motif $M$ est muni d'une filtration croissante à trois crans définie par
$$ W_{1}(M)=[0\rightarrow T]  \ \subset \ W_{2}(M)=[0\rightarrow \tilde{G}] \ \subset \  W_{3}(M)=M \: .$$
Cette filtration vérifie
$$Gr^W_{1}(M)=[0\rightarrow T], \: Gr^W_{2}(M)=[0\rightarrow A] \: \: \mathrm{et} \:\: Gr^W_{3}(M)=[Y\rightarrow 0]\: .$$ Ces gradués sont  appelés respectivement \textit{multiplicatif}, \textit{abélien}, \textit{étale}, et $W_{2}$ est appelé la partie \textit{semi-abélienne} de~$M$.

\subsubsection{Dualité, isogénie et polarisation}

Dans ce paragraphe, on définit un foncteur de dualité $M\mapsto M^t$ qui étend le foncteur $A\mapsto A^t=\Ext(A,\Gm)$ défini sur les schémas abéliens, le foncteur~$T\mapsto\Hom(T,\Gm)$ défini sur les tores et le foncteur $X\mapsto \Hom(X,\Gm)$ défini sur les réseaux.

\begin{proposition3}[\cite{HodgeIII@Deligne}, 10.2.11]  Soit $M\rightarrow S$ un $1$-motif. Le complexe
$$\tau_{\leq 0}\: \RHom_S (M,\Gm[1])$$
de $\D_S$ est concentré en degrés $-1$ et $0$ et admet un représentant canonique dans $\K_S$ qui est un $1$-motif.
\end{proposition3}

On note $M^t$ ce $1$-motif, qui dépend fonctoriellement de $M$. Comme les schémas abéliens d'une part et les tores et les $\Z$-modules d'autre part vérifient la propriété de bidualité, on peut voir que les $1$-motifs $M$ et $(M^t)^t$ sont canoniquement isomorphes. On note $A$, $T$, $\tilde{G}$, $Y$ les différents constituants de $M$, on note $X=\Hom(T,\Gm)$ et l'on note $T^t=\Hom(Y,\Gm)$. On a $$[0\rightarrow A]^t = [0\rightarrow A^t]\: ,\:\:  [0\rightarrow T]^t = [X\rightarrow 0], \:$$
$$ [Y\rightarrow 0]^t = [0\rightarrow T^t] \: \:\mathrm{et} \: \: [0\rightarrow \tilde{G}]^t = [X\buildrel{c}\over\rightarrow A^t]$$ où $c$ est l'image de l'élément de $\Ext(A,T)$ correspondant à la classe d'isomorphisme de $\tilde{G}$ par l'isomorphisme $\Ext(A,T)\iso\Hom(X,A^t)$.

Soit $\lambda  :  M_1 \rightarrow M_2$ un morphisme entre deux $1$-motifs sur $S$. On pose $\lambda_{\mathrm{m}} = Gr^W_{1}(\lambda)$, $\lambda_{\mathrm{ab}}= Gr^W_{2}(\lambda)$ et $\lambda_{\mathrm{et}} = Gr^W_{3}(\lambda)$.

\begin{definition3} On dit que $\lambda$ est une isogénie si $\lambda_{\mathrm{m}}$ et $\lambda_{\mathrm{ab}}$ sont respectivement des isogénies de tores et variétés abéliennes, et si $\lambda_{\mathrm{et}}$ est une injection à conoyau fini.
\end{definition3}

L'existence d'une isogénie entre $M_1$ et $M_2$ implique que ces derniers ont même rangs torique, abélien et étale. Notons que le dual d'une isogénie est une isogénie. Soit $\lambda : M\rightarrow M^t$ une isogénie.

\begin{definition3} On dit que $\lambda$ est une polarisation si $\lambda_{\mathrm{ab}}$ est une polarisation du gradué abélien et si $\lambda_{\mathrm{m}}^t=\lambda_{\mathrm{et}}$.
\end{definition3}

Une polarisation est \textit{principale} si c'est un isomorphisme. Par définition, le degré d'une polarisation est l'entier
$$\mathrm{deg}(\lambda)= \mathrm{deg}(\lambda_{\mathrm{ab}})\: \mathrm{deg}(\lambda_{\mathrm{m}})^2.$$
Si $M$ est polarisé, son rang torique et son rang étale sont égaux.

\subsubsection{Noyau d'une isogénie} Soit $f:M\longrightarrow M'$  une isogénie entre $1$-motifs sur $S$. On voit $f$ comme un morphisme de $\K_S$, et l'on note $\mathrm{Ker}(f)$ le $\mathrm{H}^{-1}$ du cône de
$$\xymatrix{M\ar[r]^{f} & M'}\: .$$
Il est représentable par un groupe fini et plat sur~$S$ appelé \textit{noyau} de~$f$, et est muni d'une structure d'extension panachée dont les gradués sont
\begin{eqnarray*}
Gr^W_{0}(\mathrm{Ker}(f)) &=& \mathrm{Ker}(f_\mathrm{m}) \: , \\
Gr^W_{1}(\mathrm{Ker}(f)) &=& \mathrm{Ker}(f_{\mathrm{ab}})\: , \\
Gr^W_{2}(\mathrm{Ker}(f)) &=& \mathrm{Coker}(f_{\mathrm{et}}) \: .
\end{eqnarray*}
Le dual de Cartier de~$\mathrm{Ker}(f)$ est~$\mathrm{Ker}(f^t)$.
Si $f:M\rightarrow M$ est la multiplication par un entier~$n$, on pose $M[n]=\mathrm{Ker}(f)$. Si~$T$, $A$ et $Y$ sont les constituants de~$M$, les gradués de~$M[n]$ sont $T[n]$, $A[n]$ et $Y/nY$. Comme $M[n]$ et $M^t[n]$ sont duaux de Cartier, il existe une application bilinéaire canonique $M[n]\otimes M[n]\rightarrow \mu_n$. Une polarisation principale de~$M$ induit un accouplement alterné non dégénéré $M[n]\times M^t[n] \rightarrow \mu_n$.

Soit $M$ un $1$-motif et $n$ un entier. On obtient une correspondance bijective entre l'ensemble des isogénies $f:M\rightarrow M'$ telles que $\mathrm{Ker}(f)\subset M[n]$ et l'ensemble des sous-schémas en groupes $H\subset M[n]$ en associant à~$f$ le sous-groupe $\mathrm{Ker}(f)$ et à ~$H$ l'isogénie naturelle de~$M$ vers le cône de $H[0]\rightarrow M$, qui est bien un~$1$-motif.

\subsubsection{Structures de niveau parahoriques sur les $1$-motifs}

Soient~$S$ un champ algébrique et~$M$ un~$1$-motif principalement polarisé de genre~$g$  sur~$S$. Rappelons que dans l'introduction de cet article, on a fixé un nombre premier~$p$, un entier~$s\leq g$ et un ensemble ordonné non vide
$$\mathcal{D}\:=\: \{d_1<d_2<\cdots<d_s\}\: \subset\: \{1,\cdots,g\} \: .$$

\begin{definition3} Une structure de niveau parahorique sur~$M$ de type~$\mathcal{D}$ en~$p$   est la donnée d'un drapeau  de sous-groupes finis et plats sur~$S$
$$H_1 \subset H_2 \subset \cdots \subset H_s \subset M[p]$$
tel que~$H_i$ soit de rang $p^{d_i}$ pour $i\leq s$ et que~$H_s$ soit totalement isotrope pour l'accouplement symplectique sur~$M[p]$ provenant de la polarisation.
\end{definition3}

On démontre la proposition suivante en utilisant la correspondance entre sous-groupes et isogénies et en recopiant la preuve de l'énoncé analogue pour les schémas abéliens~(\cite{Gamma0@DeJong}, proposition~1.7).

\begin{proposition3} \label{propCatPara} Les deux catégories suivantes sont équivalentes~:
\begin{itemize}
\item la catégorie des $1$-motifs~$M$ principalement polarisés de genre~$g$ munis d'une structure de niveau parahorique de type~$\mathcal{D}$ en~$p$,
\item la catégorie des diagrammes commutatifs d'isogénies entre $1$-motifs de genre $g$
$$\xymatrix{ M=M_0 \ar[d]_{} \ar[r] & M_1 \ar[r]\ar[d]_{} & \cdots \ar[r] & M_{s-1}\ar[r] \ar[d]_{} & M_s \ar[d]^{\lambda_s} \\
  M^t = M_0^t \ar[d]^{\lambda'_0} & M_1^t \ar[l] \ar[d] & \cdots \ar[l] & M_{s-1}^t \ar[d] \ar[l] & M_s^t \ar[l] \ar[d]  \\
  M=M_0 \ar[r] & M_1 \ar[r] & \cdots \ar[r] & M_{s-1} \ar[r] & M_s}$$
où la ligne horizontale médiane est duale de la ligne horizontale supérieure, la ligne horizontale inférieure est égale à la ligne horizontale supérieure, la composée $M_i\rightarrow M_i^t \rightarrow M_i$ est égale à la mutiplication par $p$ pour $i\leq s$, l'isogénie $\lambda_s$ est une polarisation, $\lambda'_0$ est l'inverse d'une polarisation principale et le degré de $M\rightarrow M_i$ est égal à $p^{d_i}$ pour $i\leq s$.
\end{itemize}
\end{proposition3}

\subsubsection{Position relative} \label{partiePosRel}

Soit $H_\bullet \subset M[p]$ une structure de niveau parahorique de type~$\mathcal{D}$. Dans ce paragraphe, on lui associe un invariant discret~$w_{H_\bullet}$ paramètrant la position relative de $H_{\bullet}$ et de la filtration de l'extension panachée $M[p]$.
Soit $V=\oplus_{j=1}^{2g}\: \F_p \cdot x_j$ un~$\F_p$-espace vectoriel de dimension~$2g$  muni d'une base~$(x_j)$. Considérons la matrice de taille~$2g \times 2g$ suivante donnée par blocs~$g\times g$:
$$J=\left( \begin{array}{cc} 0 & J' \\
 -J' & 0           \end{array} \right) \: ,$$
où~$J'$ est la matrice  anti-diagonale dont  les coefficients non nuls sont égaux à~$1$. On définit $\gsp$ comme le schéma en groupes sur $\Spec(\F_p)$ dont les~$R$-points sont
$$\gsp(R) = \{ (\gamma,\nu) \in \mathrm{GL}_{2g}(R)\times\Gm(R) \: \: | \: \:   \gamma^t\:  J\:  \gamma  \: =\:  \nu J \:\} \:.$$
Pour tout entier $r\leq g$, on définit un sous-espace totalement isotrope $W_1$ de $V$ par la formule~:
\begin{eqnarray*}
W_{1} &=& \bigoplus_{j=1}^r\:\: \F_p \cdot x_j \: .
\end{eqnarray*}
On pose $W_2=W_1^{\perp}$, $W_3=V$ et l'on obtient un drapeau à trois crans noté $W_\bullet$~;~on lui associe le sous-groupe parabolique $P=\mathrm{Stab}_{\mathrm{GSp}_{2g}}(W_\bullet)$ de $\mathrm{GSp}_{2g}$.
De même, on définit un drapeau \og standard \fg  isotrope à~$s$ crans  $V^{\mathrm{Std}_\mathcal{D}}_\bullet$ par~:
$$V^{\mathrm{Std}_\mathcal{D}}_i \: = \: \bigoplus_{j=1}^{d_i} \:\: \F_p \cdot x_j$$
pour $1\leq i \leq s$. On lui associe le sous-groupe parabolique $$P_\mathcal{D}=\mathrm{Stab}_{\mathrm{GSp}_{2g}}(V^{\mathrm{Std}_\mathcal{D}}_\bullet)$$
de $\mathrm{GSp}_{2g}$.
On pose
$$\mathcal{W}= P(\F_p) \: \backslash\ \gsp(\F_p)\: / \ P_\mathcal{D}(\F_p) \: .$$ Ce double quotient paramètre les  positions relatives de deux drapeaux de~$V$. Si~$L$ et $L_\mathcal{D}$ sont des sous-groupes de Levi de~$P$ et  $P_\mathcal{D}$, et  si $W_{\gsp}$, $W_L$ et $W_{L_\mathcal{D}}$ désignent les groupes de Weyl de $\gsp$, de $L$ et de $L_\mathcal{D}$, la décomposition de Bruhat montre que~$\mathcal{W}$ est un ensemble fini et que
$$\mathcal{W}=  W_{L} \: \backslash\ W_{\gsp} \: / \ W_{L_\mathcal{D}} \: .$$
Pour tout $w\in \mathcal{W}$ et $1\leq i \leq s$,  on pose
\begin{eqnarray*}
m_i &=& \mathrm{dim}_{\F_p} \left( (w^{}\cdot V^{\mathrm{Std}_\mathcal{D}}_i ) \: \cap \: \mathrm{Gr}_1(W_\bullet)\right) \\
a_i &=& \mathrm{dim}_{\F_p} \left(\mathrm{Im} \left( w^{}\cdot V^{\mathrm{Std}_\mathcal{D}}_i \: \cap \: W_2 \rightarrow \mathrm{Gr}_2(W_\bullet)\right)\right) \\
e_i &=& \mathrm{dim}_{\F_p} \left(\mathrm{Im} \left( w^{}\cdot V^{\mathrm{Std}_\mathcal{D}}_i \rightarrow \mathrm{Gr}_3(W_\bullet)\right)\right) \: .
\end{eqnarray*}
Ces nombres déterminent~$w$ de manière unique.

Soient~$S$ un champ algébrique connexe,~$M$ un~$1$-motif sur~$S$, principalement polarisé, de rang torique~$r$, de genre~$g$ et $H_\bullet$ une structure parahorique de type $\mathcal{D}$ en~$p$. Nous avons vu que la structure d'extension panachée définit un drapeau sur~$M[p]$, dont les gradués $\mathrm{Gr}_1^W$, $\mathrm{Gr}_2^W$ et $\mathrm{Gr}_3^W$ sont respectivement appelés multiplicatif, abélien et étale. Pour tout $1\leq i \leq s$, on pose
\begin{eqnarray*}
m_i &=& \mathrm{log}_p \: \mathrm{rg} \left(H_i \cap \mathrm{Gr}_1^W\right) \\
a_i &=& \mathrm{log}_p \: \mathrm{rg} \left(\mathrm{Im} (H_i \cap W_2(M) \rightarrow \mathrm{Gr}_2^W)\right) \\
e_i &=& \mathrm{log}_p \: \mathrm{rg} \left(\mathrm{Im} (H_i \rightarrow \mathrm{Gr}_3^W)\right)
\end{eqnarray*}
On a~$m_i+a_i+e_i=d_i$ et les suites $(m_i)_i$, $(a_i)_i$ et $(e_i)_i$ sont croissantes. On définit un drapeau isotrope à $s$ crans~$V_\bullet^{H_\bullet}$ de~$V$ en posant
$$V_i^{H_\bullet}=\left( \bigoplus_{j=1}^{m_i} \: \F_p\cdot x_j \right) \oplus \left( \bigoplus_{j=r+1}^{r+a_i} \F_p \cdot x_i \right)  \oplus  \left( \bigoplus_{j=2g-r+1}^{2g - r + e_i} \F_p\cdot x_j \right) $$
pour $1\leq i\leq s$. Par convention, si $m_i=0$ la première somme directe est nulle~;~il en est de même avec~$a_i$ et~$e_i$.
Il existe un élément (non unique) $\gamma\in \gsp(\F_p)$ qui conjugue les drapeaux $V_\bullet^{H_\bullet}$ et $V_\bullet^{\mathrm{Std}_\mathcal{D}}$. Son image dans $\mathcal{W}$ est bien déterminée et on la note $w_{H_\bullet}$~;~c'est l'invariant discret associé à la structure de niveau parahorique $H_\bullet\subset M[p]$ de type $\mathcal{D}$. Inversement, si $M\rightarrow S$ est fixé on voit aisément que pour tout $w\in \mathcal{W}$, il existe une structure de niveau $H_\bullet\subset M[p]$ telle que $w_{H_\bullet}=w$.

\subsection{Construction de Mumford}

Dans ce paragraphe, nous définissons la notion de \textit{$1$-motif de Mumford}, version plus synthétique des \textit{données de dégénérescence} du chapitre 3 de~\cite{Deg@FaltingsChai}. On note $\mathsf{Ouv}$ la catégorie dont les objets sont les diagrammes $(U\hookrightarrow S)$ où $S$ est un schéma normal et $U$ un sous-schéma ouvert dense, et dont les flèches sont les diagrammes commutatifs
$$\xymatrix{ 
U \ar[r]\ar[d] & S \ar[d] \\
U'  \ar[r] & S' \: .}$$
Soit $(U\hookrightarrow S)$ un objet de $\mathsf{Ouv}$.
\begin{definition2}\label{def1Mum} Un $1$-motif de Mumford $M$ sur $(U\hookrightarrow S)$ est un triplet $(\tilde{G},Y,M_{U})$ où
\begin{itemize} \item $\tilde{G}$ est un schéma semi-abélien de rang torique constant sur $S$, 
\item $Y$ est un réseau sur $S$,
\item $M_{U}=[Y_{U} \rightarrow \tilde{G}_{U}]$ est un $1$-motif sur $U$, où $Y_U$ et $\tilde{G}_U$ désignent les restrictions de $Y$ et de $\tilde{G}$ à $U$.
\end{itemize}
\end{definition2}
On note $\Mum_{U,S}$ la catégorie des $1$-motifs de Mumford.  La restriction à $U$ induit un foncteur $\Mum_{U,S} \rightarrow \Mot_{U}$ qui est pleinement fidèle d'après le théorème~\ref{thRay}.
Soit $\lambda \: : \: M \rightarrow M'$ un morphisme de $1$-motifs de Mumford. On dit que $\lambda$ est une \textit{isogénie} si $\lambda_U$ est une isogénie de $1$-motifs. On définit le \textit{dual} $M^t$ de $M$ comme l'unique $1$-motif de Mumford qui étend le $1$-motif $(M_U)^t$~;~son existence résulte de la construction de $M^t$. Une \textit{polarisation} $\lambda \: : \: M\rightarrow M^t$ est une isogénie telle que $\lambda_U$ soit une polarisation du $1$-motif $M_U$. On note $\Mum_{U,S,\: \mathrm{pol}}$ la catégorie des $1$-motifs de Mumford polarisés sur $(U\rightarrow S)$. Les flèches de cette catégorie sont les morphismes de $1$-motifs de Mumford qui respectent les polarisations~;~ce sont donc des isogénies.

Soit $M$ un $1$-motif de Mumford sur $(U\hookrightarrow S)$. On note $Y$, $T$, $A$, $\tilde{G}$ les schémas en groupes qui lui sont associés et $X=\Hom(T,\Gm)$. On note $\Y$ le plus grand sous-faisceau constructible de $Y$  muni d'une flèche $\Y \rightarrow \tilde{G}$ prolongeant $Y_{U} \rightarrow \tilde{G}_{U}$.  Si $\tilde{G}^t$ désigne la partie semi-abélienne de $M^t$, on note $\mathcal{X}\subset X$ le plus grand sous-faisceau constructible de $X$ muni d'une flèche $\mathcal{X}\rightarrow \tilde{G}^t$ prolongeant $X_U \rightarrow \tilde{G}^t_U$. On pose $\underline{X} = X/\mathcal{X}$ et $\underline{Y}=Y/\mathcal{Y}$~;~ce sont deux faisceaux constructibles sur $S$.
Notons $\underline{\mathrm{Div}}_S$  le faisceau des diviseurs de $S$. Nous construisons à présent une forme bilinéaire canonique
$$\underline{B} \: : \: \underline{Y} \otimes \underline{X} \longrightarrow \underline{\mathrm{Div}}_S \: .$$
Supposons d'abord que $S$ soit un trait de point générique $U=\Spec(K)$ et que les faisceaux $X$ et $Y$ soient constants. Comme $A$ est le modèle de Néron de sa fibre générique, la valuation de $K$ induit des isomorphismes
$$\tilde{G}(K)/\tilde{G}(R)\: \simeq \: T(K)/T(R)\:  \simeq \: \mathrm{Hom}(X,\Z) \: .$$
Le plongement $Y \hookrightarrow \tilde{G}(K)$ induit une forme bilinéaire $Y\otimes X\rightarrow \Z$ qui se factorise en $\underline{Y}\otimes \underline{X} \rightarrow \Z$ et ne dépend que de la classe d'isomorphisme de $M$. La construction de $\underline{B}$ dans le cas général s'en déduit par localisation le long de chaque diviseur et descente étale.

\smallskip

\begin{remarque2} La forme $\underline{B}$ prend ses valeurs dans le sous-faisceau des diviseurs de $S$ dont le support est disjoint de $U$. Si $M$ est polarisé, $\underline{B}$ vérifie la condition de symétrie suivante vis-à-vis de la partie étale $\phi : Y \rightarrow   X$ de la polarisation~:~pour toutes sections locales $y$ et $y'$ de $Y$, on a
$$\underline{B}(y,\phi(y'))=\underline{B}(y',\phi(y)) \: .$$
Si la polarisation est principale, $\underline{B}$ induit une forme bilinéaire symétrique sur $\underline{Y}$ lorsque l'on identifie $X$ et $Y$ à l'aide de $\phi$.
\end{remarque2}

\subsubsection{Construction de Mumford généralisée}

Dans cette partie, nous traduisons des résultats du chapitre~3 de~\cite{Deg@FaltingsChai} dans le langage des $1$-motifs de Mumford. Introduisons la catégorie $\mathsf{OuvComp}$ dont les objets sont les diagrammes $(U\rightarrow S \leftarrow S_0)$, où $S$ est un schéma affine normal excellent, où $S_0$ est un sous-schéma fermé défini par un idéal réduit tels que $S$ soit $S_0$-complet, et où $U$ est un sous-schéma ouvert dense de $S$ disjoint de $S_0$. Les flèches de cette catégorie sont les diagrammes commutatifs
$$\xymatrix{ 
U \ar[r]\ar[d] & S \ar[d] & S_0 \ar[l] \ar[d] \\
U'  \ar[r] & S' & S'_0 \ar[l]}$$

Introduisons maintenant une notion de positivité sur les $1$-motifs de Mumford polarisés. Soient $M$ un $1$-motif de Mumford polarisé sur $(U\hookrightarrow S)$ et $\phi  :  Y \rightarrow X$ l'injection à conoyau fini associée.

\begin{definition3} \label{positivite} $M$ est positif relativement à $S_0$ si pour toute section locale non nulle $y$ de $Y$, le diviseur $\underline{B}(y,\phi(y))$ est effectif et son support contient $S_0$.
\end{definition3}

La positivité de $M$ implique que $\mathcal{Y}_{S_0}=\mathcal{X}_{S_0}=0$ et que $M_{U}$ est un $1$-motif cohomologiquement concentré en degré $0$, \textit{i.e.} $H^{-1}(M_{U})=0$ et $Y_{U}$ s'injecte dans $\tilde{G}_{U}$. On note $$\Mum_{U,S,S_0,\: pol,\: +}$$ la catégorie des $1$-motifs de Mumford polarisés sur $(U\hookrightarrow S)$ positifs relativement à $S_0$. La proposition~\ref{propDel} montrera qu'elle est équivalente à la catégorie $DD_{pol}$ définie dans le paragraphe~III.2 de~\cite{Deg@FaltingsChai}. 
On note $\mathsf{Deg}_{U,S,S_0,\: pol}$ la catégorie des schémas semi-abéliens sur $S$, abéliens polarisés sur $U$ et de rang torique constant sur $S_0$. On demande en outre que la partie torique de $G_{S_0}$ soit isotriviale sur $S_0$. Les flèches sont les morphismes de schémas semi-abéliens qui respectent la polarisation~;~ce sont donc des isogénies.

\begin{remarque3} On ne trouve pas la condition d'isotrivialité dans~\cite[ch.~II]{Deg@FaltingsChai}. Comme le remarque Lan~\cite[3.3.3.1]{Compact@Lan}, elle est pourtant nécessaire pour algébriser l'extension de Raynaud de $G$ (\textit{cf.} aussi~\sga{7}{ix}{7.2.1} et~\sga{3}{x}{3.2}). De plus, cette condition permet d'éliminer l'hypothèse~\cite[II.3.a]{Deg@FaltingsChai}.
\end{remarque3}

Le but de la construction de Mumford est d'exhiber un foncteur $$\mathrm{M}_{pol} \: : \: \Mum_{U,S,S_0,\:pol, \: +} \longrightarrow \mathsf{Deg}_{U,S,S_0,\: pol}$$ tel que si $S$ est trait complet de point fermé $S_0$  et de point générique $U=\Spec(K)$, on ait ensemblistement
$$\mathrm{M}_{pol}(M)(K)=H^0\left(M(K)\right).$$

\begin{theoreme3}[\cite{Deg@FaltingsChai}, corollaires III.7.2, 7.3 et 7.5] \label{ThConsMum}  Il existe un foncteur canonique explicite
$$\mathrm{M}_{pol} \: : \: \Mum_{U,S,S_0,\: pol, \: +} \: \longrightarrow \: \mathsf{Deg}_{U,S,S_0,\: pol} \: .$$
Il réalise une équivalence de catégories. Soit $M=(\tilde{G},Y,M_U=[Y_U \rightarrow \tilde{G}_U])$. Si l'on note $G=\mathrm{M}_{pol}(M)$, et si l'on désigne par $\underline{X}(G)$ le faisceau constructible associé à $G$ par le théorème~\ref{thFaisCons}, et par $\underline{X}$ le faisceau constructible associé à $M$, on  dispose d'isomorphismes fonctoriels canoniques
\begin{eqnarray*}  G_{S_{0}} & \isolong&   \tilde{G}_{S_{0}}, \\
\underline{X}(G) & \isolong &\underline{X}, \\
G_{U}[p] & \isolong & M_{U}[p] \: .
\end{eqnarray*}
En particulier, $G_U[p]$ est muni d'une structure d'extension panachée.
\end{theoreme3}

\subsubsection{Globalisation de la forme bilinéaire symétrique} \label{parGlobB}

Soient $S$ un champ algébrique normal excellent, $U$ un ouvert dense de $S$ et $G$ un schéma semi-abélien sur $S$  abélien polarisé sur $U$. On note $K_{U}$ le noyau de la polarisation de $G_{U}$ et $K$ son adhérence dans $G$. Cette adhérence est  un schéma en groupes quasi-fini et plat sur $S$. Dans la partie II.2 de~\cite{Deg@FaltingsChai}, Faltings et Chai montrent que le quotient $G/K$ est représentable par un schéma semi-abélien  noté $G^t\rightarrow S$. Comme $G_U^t$ est indépendant du choix de la polarisation, le théorème~\ref{thRay} montre que $G^t$ ne dépend pas du choix de la polarisation de $G_U$.

Pour tout objet $(U\rightarrow S \leftarrow S_0)$ de $\mathsf{OuvComp}$, la construction de Mumford est compatible à la dualité dans le sens où si $M$ est objet de $\Mum_{U,S,S_0,\:pol, \: +}$ tel que $G=\mathrm{M}_{pol}(M)$, on a $G^t\iso\mathrm{M}_{pol}(M^t)$. Dans la proposition suivante, on se place à nouveau dans le cas où $S$ est un champ algébrique normal excellent quelconque et $U$ est un ouvert de $S$.

\begin{proposition3}[\cite{Deg@FaltingsChai}, chapitre III, théorème 10.1] \label{propBGsab}  Il existe une unique forme bilinéaire $$\underline{B}(G) \: : \: \underline{X}(G^t) \times \underline{X}(G) \longrightarrow \underline{\mathrm{Div}}_S$$ telle qu'après tout changement de base par un anneau local intègre complet $S'\rightarrow S$ on ait
$$\underline{B}(G)=\mathrm{M}_{pol}\left(\underline{B}(\mathrm{M}_{pol}^{-1}(G)\right) \: .$$
\end{proposition3}

D'après le théorème~\ref{thRay}, la polarisation $G_\eta \rightarrow G_\eta ^t$ s'étend en un morphisme $G\rightarrow G^t$ qui induit une flèche $\phi \: :\: \underline{X}(G^t) \rightarrow \underline{X}(G)$. La forme bilinéaire $\underline{B}$ vérifie la relation
$$\underline{B}\left(y,\phi(y')\right)=\underline{B}\left(y',\phi(y)\right)$$
pour toutes sections locales $y$ et $y'$ de $\underline{X}(G^t)$ et la positivité de $\mathrm{M}_{pol}^{-1}(G)$ implique que $\underline{B}(y,\phi(y))$ est un diviseur effectif dont le support est le lieu de non annulation de $y$.

\subsubsection{Application de Kodaira-Spencer}

Soient $(U\rightarrow S \leftarrow S_0)$  un objet de $\mathsf{OuvComp}$ et~$G$ un objet de $\mathsf{Deg}_{U,S,S_0,\: pol}$. On note $G^t$ le dual de $G$ défini dans le paragraphe précédent, $\lambda  :  G \rightarrow G^t$ l'extension sur $S$ de la polarisation sur $U$, on note $e_G :  S \rightarrow G$ et $e_{G^t}  : S \rightarrow G^t$ les sections neutres, on note $\Omega_G=e_G^* \: \Omega_{G/S}$ le faisceau des $1$-formes différentielles invariantes sur $G$ et l'on note $\Omega_{G^t}=e_{G^t}^* \: \Omega_{G^t/S}$. On suppose que $S\rightarrow \Spec(\Z)$ est formellement lisse et que le complémentaire $\partial S$ de $U$ dans $S$ est un diviseur à croisements normaux.

\begin{proposition3}[\cite{Deg@FaltingsChai}, partie III.9, page 84]  \label{propKSMum} Il existe une flèche canonique
$$K_{M/S} \:~:~\: \Omega_G \otimes \Omega_G^t \longrightarrow \Omega^1_{S/\Z} (\mathrm{log}\: \partial  S)$$appelée flèche de Kodaira-Spencer. Elle se factorise en
$$K_{M/S} \: : \: \mathrm{Coker}\left(\Omega_{G^t} \otimes \Omega_{G^t} \longrightarrow \Omega_G \otimes \Omega_G^t \right) \longrightarrow \Omega^1_{S/\Z} (\mathrm{log}\: \partial  S) \: ,$$
où la flèche $\Omega_{G^t} \otimes \Omega_{G^t} \longrightarrow \Omega_G \otimes \Omega_G^t$ envoie $x\otimes y$ sur $\lambda^*(x) \otimes y - \lambda^*(y)\otimes x$.
\end{proposition3}

La flèche de Kodaira-Spencer est construite grâce aux $1$-motifs obtenus par l'inverse de la construction de Mumford. D'après~\cite[coro.~III.9.8]{Deg@FaltingsChai}, elle se globalise dans le cas où $(U\hookrightarrow S)$ est un objet de $\mathsf{Ouv}$ et $G$ un schéma semi-abélien sur $S$ qui est abélien polarisé sur~$U$.

\subsubsection{Structure de niveau parahorique}

Si $(U\rightarrow S \leftarrow S_0)$ est un objet de $\mathsf{OuvComp}$, on note $\mathsf{Deg}_{U,S,S_0,\: pol, \: 0}$ la catégorie des schémas semi-abéliens $G$ sur $S$ de partie torique isotriviale sur $S_0$  dont la restriction à $U$ est un schéma abélien  polarisé muni d'une structure de niveau parahorique de type $\mathcal{D}$ en $p$. De même, on note $\Mum_{U,S,S_0,\: pol, \: +, \: 0}$ la catégorie des $1$-motifs de Mumford polarisés positifs dont le $1$-motif sous-jacent sur $U$ est muni d'une structure de niveau parahorique de type $\mathcal{D}$ en $p$.
La proposition suivante découle du théorème \ref{ThConsMum} par fonctorialité de la construction de Mumford.

\begin{proposition3} \label{propMumIw} Le foncteur $\mathrm{M}_{pol}$ induit une équivalence de catégories $$\mathrm{M}_{pol}\: : \: \Mum_{U,S,S_0,\: pol, \: +, \: 0} \: \longrightarrow \: \mathsf{Deg}_{U,S,S_0,\: pol, \: 0}$$
\end{proposition3}

Cette proposition justifie l'étude des espaces de modules de $1$-motifs polarisés munis de structures de niveau parahorique.

\subsection{Espace de modules de $1$-motifs}

Dans ce paragraphe, on montre que les champs de modules de $1$-motifs principalement polarisés qui sont munis de structures de niveau parahorique sont algébriques, et on les décrit explicitement. Commençons par traiter le cas d'une structure de niveau triviale.

\subsubsection{Description symétrique d'un $1$-motif}

Pour tout schéma abélien $A$, on note $\mathcal{P}_A$ le fibré de Poincaré sur $A\times A^t$.

\begin{proposition3}[\cite{HodgeIII@Deligne}, proposition~10.2.14] \label{propDel} La catégorie des $1$-motifs est équivalente à la catégorie dont les objets sont les familles formées de 
\begin{itemize} \item un schéma abélien $A$,
\item deux faisceaux isotriviaux $X$ et $Y$,
\item des morphismes $c  : X \rightarrow A^t$ et $c^t :  Y \rightarrow A$,
\item une trivialisation de la $\Gm$-biextension $(c^t \times c)^{*} \mathcal{P}_A$ sur $Y \times X$.
\end{itemize}
\noindent
La catégorie des $1$-motifs principalement polarisés est équivalente à la catégorie dont les objets sont les familles formées de
\begin{itemize} \item un schéma abélien $A$ principalement polarisé par $\lambda_A$,
\item un réseau $X$,
\item un morphisme $c : X \rightarrow A$,
\item une trivialisation symétrique de la $\Gm$-biextension symétrique $(c\times c)^* \mathcal{P}_A$, où  $\mathcal{P}_A$ est vu comme une $\Gm$-biextension sur $A\times A$ \textit{via} $\lambda_A$.
\end{itemize}
\end{proposition3}

\subsubsection{Espace de modules}

Soient $r\leq g$ un entier et $\mathbb{X}$ un $\Z$-module libre de rang $r$. On note $\mathbb{T}=\Hom(\mathbb{X},\Gm)$, on note $\mathcal{A}_{g-r}$ le champ des schémas abéliens principalement polarisés de genre $g-r$, et l'on note $\mathrm{A}_{g-r} \rightarrow \mathcal{A}_{g-r}$ le schéma abélien universel.

\begin{definition2} Un $1$-motif principalement polarisé rigidifié par $\mathbb{X}$ est la donnée d'une famille $(M,\lambda,\mathrm{Gr}^W_1(M) \iso \mathbb{X},\mathrm{Gr}^W_3(M) \iso \mathbb{T})$ où $M$ est un $1$-motif, $\lambda$ une polarisation principale de $M$ et la composée de  $\mathrm{Gr}^W_1(\lambda)$ avec les isomorphismes $\mathrm{Gr}^W_1(M) \iso \mathbb{X}$ et $\mathrm{Gr}^W_1(M^t)=\Hom(\mathrm{Gr}^W_3(M),\Gm)\iso \mathbb{X}$ est l'identité de $\mathbb{X}$.
\end{definition2}

Comme le groupe d'automorphismes de la partie abélienne polarisée $\mathrm{Gr}^W_2(M)$
est fini, le groupe d'automorphismes d'un $1$-motif principalement polarisé rigidifié est également fini. On note $\mathcal{M}$ le champ des $1$-motifs principalement polarisés de genre $g$ qui sont rigidifiés par $\mathbb{X}$.

\begin{proposition2} \label{propB} Le champ $\mathcal{M}$ est algébrique et se dévisse en
$$\mathcal{M}\longrightarrow  \mathcal{B} = \Hom(\mathbb{X},\mathrm{A}_{g-r}) \longrightarrow \Spec(\Z)\: ,$$
où $\mathcal{M}\rightarrow \mathcal{B}$ est un torseur sous le tore $E=\Hom\left(\mathrm{Sym}^2(\mathbb{X}),\Gm\right)$.
\end{proposition2}
La classe d'isomorphisme de ce torseur est donnée dans la démonstration suivante.
\medskip

\begin{demo} Il suffit de montrer la deuxième assertion. On utilise la proposition~\ref{propDel} en réutilisant les mêmes notations. La flèche $\mathcal{M}\rightarrow \mathcal{B}$ envoie $M$ sur $c  :  \mathbb{X} \rightarrow \mathrm{A}_{g-r}$ et $\mathcal{M}\rightarrow \mathcal{B}$ est l'espace de modules relatif des trivialisations de la biextension symétrique $(c\times c)^* \mathcal{P}_{\mathrm{A}_{g-r}}$. La classe d'isomorphisme de cette biextension définit un élément de $$\Ext_{\mathcal{B}}\left(\mathrm{Sym}^2(\mathbb{X}),\Gm\right)\: =\: \Ext_{\mathcal{B}}\left(\Z,\Hom\left(\mathrm{Sym}^2(\mathbb{X}),\Gm\right)\right)$$
donc une classe d'isomorphisme de  $E$-torseurs sur $\mathcal{B}$, qui est celle de $\mathcal{M}$.
\end{demo}

\subsubsection{Description symétrique d'une structure de niveau}

Donnons une variante de la proposition~\ref{propDel} pour les $1$-motifs munis d'une structure de niveau parahorique de type $\mathcal{D}$ en $p$. Soit $w\in \mathcal{W}$ d'entiers associés $m_i$, $a_i$ et $e_i$.

\begin{proposition3} \label{propParaSym} Le champ des $1$-motifs principalement polarisés de genre $g$, de rang torique $r$,  munis d'une structure de niveau parahorique de type $\mathcal{D}$ et d'invariant $w$ en $p$ est équivalent au champ qui paramètre les familles formées d'un diagramme commutatif
$$\xymatrix{
&&&& Y_0\ar[lldd]\ar[d]\ar[rr] & & Y_1\ar[lldd]\ar[d] \ar[rr] & & \cdots \ar[d] \ar[lldd] \ar[rr] & & Y_s \ar[d] \ar[lldd] \\
&&&& A_0\ar'[ld][lldd]\ar'[r][rr] & & A_1\ar'[ld][lldd] \ar'[r][rr] & & \cdots \ar'[ld][lldd] \ar'[r][rr] & & A_s \ar[lldd] \\
&& X_0\ar[lldd] \ar[d] &  & X_1\ar[ll] \ar[d]\ar[lldd]  & & \cdots \ar[d] \ar[lldd] \ar[ll] & & X_s \ar[d] \ar[ll] \ar[lldd]\\
&& A_0^t \ar'[ld][lldd] & & A_1^t\ar'[l][ll]\ar'[ld][lldd]   && \cdots \ar'[ld][lldd] \ar'[l][ll] & & A_s ^t \ar[lldd] \ar'[l][ll] \\
Y_0 \ar[rr]\ar[d] & & Y_1\ar[d]\ar[rr] & & \cdots  \ar[d] \ar[rr]  & & Y_s \ar[d] \\
A_0\ar[rr] & & A_1\ar[rr] && \cdots \ar[rr] & & A_s}$$ 
et de trivialisations $\tau_i$ de la $\Gm$-biextension $(Y_i\times X_i \rightarrow A_i\times A_i^t)^*\: \mathcal{P}_i$ sur $Y_i\times X_i$ pour  $0\leq i\leq s$ (où $\mathcal{P}_i$ est la biextension de Poincaré sur $A_i\times A_i^t$)
telles que~:
\begin{itemize} 
\item la face avant et la face arrière du diagramme sont égales,
\item les objets de la face supérieure sont des faisceaux isotriviaux et les flèches de cette face sont des injections~;~$X_0\rightarrow Y_0$ est un isomorphisme~;~pour tout $i\leq s$, les composés $Y_i\rightarrow X_i\rightarrow X_i$ sont égaux à la multiplication par $p$ et les conoyaux des composés $Y_0\rightarrow Y_i$ et $X_i\rightarrow X_0$ sont de cardinaux respectifs $p^{e_i}$ et $p^{m_i}$,
\item les objets de la face inférieure sont des schémas abéliens de genre $g-r$ et les flèches de cette face sont des isogénies~;~le complexe $A_0^t \leftarrow A_1^t \leftarrow \cdots \leftarrow A_s^t$ est le dual du complexe $A_0\rightarrow A_1\rightarrow \cdots \rightarrow A_s$~;~$A_0\rightarrow A_0^t$ et $A_s\rightarrow A_s^t$ sont des polarisations, $A_0^t \rightarrow A_0$ est un isomorphisme~;~pour tout $i\leq s$, la composée $A_i\rightarrow A_i^t \rightarrow A_i$ est la multiplication par $p$ et le degré de l'isogénie  $A_0\rightarrow A_i$ est $p^{a_i}$.
\item pour tout $i\leq s$, la trivialisation $\tau_i$ induit une trivialisation symétrique de la biextension symétrique $(Y_i\times Y_i \rightarrow Y_i \times X_i \rightarrow A_i\times A_i^t)^* \: \mathcal{P}_i$,
\item les trivialisations $\tau_i$ sont compatibles aux flèches de transitions $Y_i\rightarrow Y_{i+1}$, $X_i \leftarrow X_{i+1}$ et $A_i\rightarrow A_{i+1}$.
\end{itemize}
\end{proposition3}

\subsubsection{Espace de modules avec structure de niveau parahorique} \label{partieEspModPar}
 
On pose $\mathbb{Y}=\mathbb{X}$ et l'on note $\mathcal{M}_{0}$ le champ des $1$-motifs de genre $g$, principalement polarisés et rigidifiés par $\mathbb{Y}$ et munis d'une structure parahorique de type $\mathcal{D}$ en $p$. On va montrer qu'il est algébrique et en donner une description explicite. On peut décomposer $\mathcal{M}_0$ en l'union disjointe
$$\mathcal{M}_0 = \coprod_{w\in\mathcal{W}}  \mathcal{M}^w_0$$
où $\mathcal{M}^w_0$ est le champ des $1$-motifs munis d'une structure de niveau parahorique d'invariant $w$, qui est ouvert et fermé dans $\mathcal{M}_0$. Il suffit de montrer que $\mathcal{M}^w_0$ est algébrique pour tout $w\in \mathcal{W}$, et de le décrire explicitement.
On note $\mathcal{B}_0^w$ le champ de modules des diagrammes commutatifs comme dans la proposition~\ref{propParaSym}, qui sont munis d'un isomorphisme entre 
$$\xymatrix{Y_0\ar[r] &  X_0 \ar[r] & Y_0 } \:\: \mathrm{et} \:\: \xymatrix{\mathbb{Y} \ar[r]^{\times \: p} & \mathbb{X} \ar[r]^{\mathrm{Id}} &\mathbb{Y}}\: .$$
On remarque que les degrés des flèches $Y_i\rightarrow Y_{i+1}$, $X_{i+1}\rightarrow X_i$ et $A_i\rightarrow A_{i+1}$ sont déterminés par $w$.

\begin{proposition3} Le champ algébrique $\mathcal{B}_0^w$ est un schéma en groupes propre sur une union disjointe de variétés de Siegel $\agrK$ dont le type parahorique dépend de $w$.
\end{proposition3}

\begin{demo} Soit $\mathcal{B}'^w_0$ le champ de modules des couples de diagrammes commutatifs ayant les mêmes propriétés que le couple formé par la face supérieure et la face inférieure du diagramme de la proposition~\ref{propParaSym}. Il est algébrique et est une union disjointe de variétés de Siegel $\agrK$ de genre $g-r$ 
et de type parahorique $$a_1\leq a_2 \leq \cdots \leq a_g \leq g-r \: .$$
Le champ $\mathcal{B}^w_0\rightarrow \mathcal{B}'^w_0$ paramètre les flèches $Y_i\rightarrow A_i$ et $X_i\rightarrow A_i^t$ qui vérifient des conditions de commutation adéquates. Ce champ est égal à la limite projective du diagramme suivant~:

\newpage

{\small \begin{equation} \label{equationLosange} \xymatrix@C=0pt{
& \Hom_{\mathcal{B}'^w_0}(X_0,A_0) & \\
\Hom_{\mathcal{B}'^w_0}(X_0,A_0^t) \ar[ru] \ar[rd] \ar[dd] & & \Hom_{\mathcal{B}'^w_0}(Y_0,A_0) \ar[ld] \ar[lu] \ar[dd]  \\
 & \Hom_{\mathcal{B}'^w_0}(Y_0,A_0^t ) & \\
\Hom_{\mathcal{B}'^w_0}(X_1,A_0^t) & & \Hom_{\mathcal{B}'^w_0}(Y_0,A_1)  \\
& \Hom_{\mathcal{B}'^w_0}(X_1,A_1)  & \\
\Hom_{\mathcal{B}'^w_0}(X_1,A_1^t) \ar[dd] \ar[uu] \ar[ur] \ar[dr] & & \Hom_{\mathcal{B}'^w_0}(Y_1,A_1) \ar[uu] \ar[lu] \ar[ld] \ar[dd]\\
& \Hom_{\mathcal{B}'^w_0}(Y_1,A_1^t) & \\
\cdots   & & \cdots  \\
& \Hom_{\mathcal{B}'^w_0}(X_s,A_s) & \\
\Hom_{\mathcal{B}'^w_0}(X_s,A_s^t) \ar[ru] \ar[rd]\ar[uu]  & & \Hom_{\mathcal{B}'^w_0}(Y_s,A_s) \ar[ld] \ar[lu]  \ar[uu] \\
 & \Hom_{\mathcal{B}'^w_0}(Y_s,A_s^t ) & }
\end{equation}}

Cela montre que le morphisme $\mathcal{B}^w_0\rightarrow \mathcal{B}'^w_0$ est représentable par un schéma en groupes propre et que le champ $\mathcal{B}^w_0$ est algébrique.
\end{demo}

La proposition~\ref{propParaSym} implique que le champ relatif $\mathcal{M}^w_0\rightarrow \mathcal{B}^w_0$ paramètre certaines familles compatibles de trivialisations de $\Gm$-biextensions. On  en déduit le résultat suivant.

\begin{proposition3} \label{propXi0w} La flèche $\mathcal{M}^w_0\rightarrow \mathcal{B}^w_0$ est un torseur sous un groupe de type multiplicatif $E_0^w$ isogène à $E$. En particulier, le champ $\mathcal{M}_0^w$ est algébrique.
\end{proposition3}

L'analyse des deux dernières conditions de la proposition~\ref{propParaSym} montre que le groupe des caractères de $E_0^w$ est la limite inductive du diagramme suivant~:
\begin{equation} \label{equationFlip} \xymatrix{
X_0 \otimes X_0  \ar@<1ex>[d]|f \ar@<-1ex>[d] &  & X_1\otimes X_1 \ar@<-1ex>[d] \ar@<1ex>[d]|f  & &  X_s\otimes X_s \ar@<-1ex>[d] \ar@<1ex>[d]|f \\
Y_0\otimes X_0 & Y_0\otimes X_1 \ar[r]\ar[l] & Y_1\otimes X_1  & \cdots \ar[r] \ar[l] & Y_s\otimes X_s \\
Y_0\otimes Y_0 \ar@<1ex>[u]|f \ar@<-1ex>[u] &  & Y_1\otimes Y_1 \ar@<1ex>[u]|f \ar@<-1ex>[u] & & Y_s\otimes Y_s \ar@<1ex>[u]|f \ar@<-1ex>[u]}
\end{equation}
Les flèches sont induites par celles du diagramme
\begin{equation} \label{equationXY} \xymatrix{ Y_0 \ar[r] \ar[d] & Y_1 \ar[r] \ar[d] & \cdots  \ar[r] & Y_s \ar[d] \\
X_0 \ar[d]  & X_1 \ar[l]\ar[d] & \cdots \ar[l] & X_s\ar[d] \ar[l] \\
Y_0 \ar[r] & Y_1 \ar[r] & \cdots \ar[r]  & Y_s }
\end{equation}
et $\xymatrix{\ar[r]|f &}$ signifie que l'on a composé la flèche correspondante avec l'isomorphisme de~\emph{flip} (ou symétrie) qui envoie $a\otimes b$ sur $b \otimes a$.

\subsubsection{Liberté} Dans ce paragraphe, nous montrons que  $\mathcal{B}^w_0\rightarrow \mathcal{B}'^w_0$ est un schéma abélien et que $E_0^w$ est un tore. Ces résultats seront utilisés  dans le corollaire~\ref{coroXsiKbNormal}, lui-même nécessaire pour démontrer la proposition~\ref{propR0bet}. 
On se donne des symboles $\varepsilon_k$ pour $1\leq k \leq r$ et on note $[k,k']$ l'ensemble des entiers compris entre $k$ et $k'$. 
Soient $e_i$, $m_i$ et $a_i$ les entiers associés à $w\in \mathcal{W}$. On fixe  des isomorphismes $Y_i\iso \oplus_{k=1}^r \Z \cdot \varepsilon_k$ et $X_i\iso \oplus_{k=1}^r \Z \cdot \varepsilon_k$ pour $0\leq i \leq s$, tels que
\begin{itemize}
\item le morphisme $Y_0 \rightarrow Y_i$ envoie $\varepsilon_k$ sur $p\: \varepsilon_k$ si $k\in [1,e_i]$, et sur  $\varepsilon_k$ sinon,
\item le morphisme $X_s \rightarrow X_i$ envoie $\varepsilon_k$ sur $p\: \varepsilon_k$ si $k\in [r-m_s+1,r-m_i]$, et sur  $\varepsilon_k$ sinon.
\end{itemize}

\paragraph{Cas de $\mathcal{B}^w_0$} Commençons par énoncer quelques lemmes préliminaires.

\begin{lemme3}\label{lemmeAbFormel} Soient $B$ un schéma abélien, $H_1$, $H_2$, $H$ et $H'$ quatre sous-groupes finis et plats de $B$ tels que $H_1+H_2 \: \subset \: H\cap H'$. La limite projective  de
$$\xymatrix{ B/H_1 \ar[r]\ar[d] & B/H' \\
B/H & B/H_2\ar[l]\ar[u] }$$
est extension de $H\cap H' \: / \: H_1 + H_2$ par $B\:/\: H_1\cap H_2$.
\end{lemme3}

\begin{demo} Cette limite projective $P$ est égale à $$\lbrace (b_1,b_2)\in B\times B \: | \: b_1-b_2 \in H \cap H'  \rbrace \: / \: H_1\times H_2\: .$$ La flèche qui envoie $(b_1,b_2)$ sur $b_1-b_2$ définit une surjection de $P$ sur $H\cap H' \: / \: H_1 + H_2$ , dont le noyau est
$$ \lbrace (b_1,b_2)\in B\times B \: | \: b_1-b_2 \in H_1+H_2  \rbrace \: / \: H_1\times H_2\: .$$
On conclut en remarquant que l'immersion diagonale de $B$ dans $B\times B$ induit un isomorphisme de $B\:/\: H_1\cap H_2$ sur ce noyau.
\end{demo}

\begin{lemme3} \label{lemmeAbLosange} La limite projective $P'_i$ du diagramme
$$\xymatrix{\Hom(Y_i,A_i) \ar[r]\ar[d] & \Hom(X_i,A_i)  \\
\Hom(Y_i,A_i^t) & \Hom(X_i,A_i^t) \ar[u] \ar[l]}$$
est représentable par un schéma abélien, qui est égal au quotient de $\Hom(Y_i,A_i)$ par le sous-groupe fini et plat formé des éléments $\phi$ de $\Hom(Y_i,A_i)$ tels que $\phi (\varepsilon_k) \in \mathrm{Ker}(A_i\rightarrow A_i^t)$ si $k\in [e_i +1 , r-m_i]$ et $\phi(\varepsilon_k) \in A_i[p]$  sinon.
\end{lemme3}
 
\begin{demo} On pose $B = \Hom(Y_i,A_i)$ et
$$\begin{array}{lcl}
H_1 & = &  \lbrace \:\phi \in B \:\: | \:\: \phi(p Y_i) \:=\: 0 \:\rbrace, \\
H_2 & = & \lbrace \:\phi \in B \:\: | \:\:  \phi(X_i) \:\subset \:\mathrm{Ker}(A_i\rightarrow A_i^t)\:\rbrace\: ,\\
H & = & \lbrace\: \phi \in B  \: \:| \: \: \phi(p Y_i)\: \subset \:\mathrm{Ker}(A_i\rightarrow A_i^t) \:\rbrace \: , \\
H' & = & \lbrace \:\phi \in B \: \: | \: \: \phi(X_i) \:\subset \:A_i[p] \:\rbrace\:.
\end{array}$$
On applique le lemme~\ref{lemmeAbFormel} puis l'on observe que
$H\cap H' = H_1 + H_2$ et que $H_1\cap H_2$ est le schéma en groupes annoncé, qui est bien sûr fini. On voit qu'il est plat grâce à sa description explicite.
\end{demo}

Pour tout $0\leq i\leq s$, on note $P_i$ la limite projective du diagramme obtenu à partir de~(\ref{equationLosange}) en ne conservant que les objets dont l'indice est~$\leq i$.

\begin{lemme3} \label{lemmeAbRec} Le morphisme naturel de $\Hom(Y_i,A_0)$ dans $P_i$ est surjectif. Son noyau est formé des éléments $\phi \in \Hom(Y_i,A_0)$ vérifiant 
$\phi(\varepsilon_k)\in\mathrm{Ker}(A_0\rightarrow A_i^t)$ {pour} $k\in[e_i+1,r-m_i]$.
En particulier, $P_i$ est un schéma abélien.
\end{lemme3}

\begin{demo} On procède par récurrence sur $i$. Le cas $i=0$ résulte du lemme~\ref{lemmeAbLosange}. Supposons que les assertions concernant $P_{i-1}$ soient vérifiées et prouvons celles concernant $P_i$. Le groupe $P_i$ est la limite projective de
$$\xymatrix{P_{i-1}\ar[r]\ar[d] & \Hom(Y_{i-1},A_i) \\
\Hom(X_i,A_{i-1}^t) & P'_i\ar[l]\ar[u] }$$
Il est quotient de la limite projective de
$$\xymatrix{\Hom(Y_{i-1},A_0)\ar[r]\ar[d] & \Hom(Y_{i-1},A_i) \\
\Hom(X_i,A_{i-1}^t) & \Hom(Y_i,A_i)\ar[l]\ar[u] }$$
par un sous-groupe dont la projection $Q_i$ sur $\Hom(Y_{i-1},A_0)$ est formée de morphismes vérifiant 
$\phi(\varepsilon_k)\in\mathrm{Ker}(A_0\rightarrow A_{i-1}^t)$ pour $k\in[e_{i-1}+1,r-m_{i-1}]$. Posons $B=\Hom(Y_i,A_0)$ et
$$\begin{array}{lcl}
H_1 & = &  \lbrace \:\phi \in B \:\: | \:\: \phi( Y_{i-1}) \:= \:0 \:\rbrace, \\
H_2 & = & \lbrace \:\phi \in B \:\: | \:\:  \phi(Y_i) \:\subset \:\mathrm{Ker}(A_0\rightarrow A_i)\:\rbrace\: ,\\
H & = & \lbrace \:\phi \in B  \:\: | \: \: \phi( Y_{i-1}) \:\subset \:\mathrm{Ker}(A_0\rightarrow A_i)\:  \rbrace \: , \\
H' & = & \lbrace\: \phi \in B \: \: | \: \: \phi(X_i) \:\subset \: \mathrm{Ker}(A_0\rightarrow A_{i-1}^t)  \:\rbrace\:.
\end{array}$$
Si $R_i$ désigne l'image inverse de $Q_i$ dans $B$, on voit que $H\cap H' \: / \: (H_1 + H_2 + R_i)=0$. Le lemme~\ref{lemmeAbFormel} montre alors la surjectivité du morphisme naturel de $B$ vers $P_i$. L'assertion concernant son noyau résulte du calcul de $H_1\cap H_2$, ce qui conclut la récurrence.
\end{demo}

\begin{corollaire3} \label{coroAb} Le schéma  $\mathcal{B}^w_0 $ est un schéma abélien sur $\mathcal{B}'^w_0$.
\end{corollaire3}

\paragraph{Cas de $E_0^w$} Soit $S_0^w$ le groupe des caractères de $E_0^w$. Il est limite inductive du diagramme~(\ref{equationFlip}).

\begin{lemme3} \label{lemmeEtCran} La limite inductive $S'_i$ du diagramme
$$\xymatrix{ X_i \otimes X_i  \: \ar@<0.5ex>[r]|f \ar@<-0.5ex>[r] &\: Y_i\otimes X_i \: & \: Y_i \otimes Y_i \ar@<0.5ex>[l]|f \ar@<-0.5ex>[l]}$$
est un $\Z$-module libre. Les éléments $\varepsilon_k\otimes \varepsilon_l $ et  $\varepsilon_l\otimes \varepsilon_k \in Y_i\otimes X_i$ ont même image dans $S'_i$ pour $(k,l)\in [e_i +1,r-m_i]^2$.
\end{lemme3}

\begin{demo} Le morphisme $Y_i\otimes Y_i \rightarrow Y_i\otimes X_i$ obtenu par différence envoie $\varepsilon_k\otimes \varepsilon_l$ sur 
$$\begin{array}{lcl}
\varepsilon_k\otimes\varepsilon_l - \varepsilon_l \otimes \varepsilon_k   & \mathrm{si} &  (k,l)\in[e_i+1,r-m_i]^2 \: ,\\
p \: \varepsilon_k\otimes\varepsilon_l - \varepsilon_l \otimes \varepsilon_k &  \mathrm{si} &   (k,l)\in \left([1,e_i]\cup[r-m_i+1,r]\right)\times [e_i+1,r-m_i]\: , \\
\varepsilon_k\otimes\varepsilon_l - p\: \varepsilon_l \otimes \varepsilon_k  & \mathrm{si} &   (k,l)\in [e_i+1,r-m_i]\times \left([1,e_i]\cup[r-m_i+1,r]\right) \: ,\\
p\: \varepsilon_k\otimes\varepsilon_l - p \: \varepsilon_l \otimes \varepsilon_k & \mathrm{si} &   (k,l)\in \left([1,e_i]\cup[r-m_i+1,r]\right) \times \left([1,e_i]\cup[r-m_i+1,r]\right) .
\end{array}$$
D'autre part, le morphisme $X_i\otimes X_i \rightarrow Y_i\otimes X_i$ envoie $\varepsilon_k\otimes \varepsilon_l$ sur $\varepsilon_k\otimes\varepsilon_l - \varepsilon_l\otimes\varepsilon_k$ si $k$ et $l$ sont dans  $\left([1,e_i]\cup[r-m_i+1,r]\right)$. Cela suffit pour conclure.
\end{demo}

Pour tout $0\leq i\leq s$, on note $S_i$ la limite inductive du diagramme obtenu à partir de~(\ref{equationFlip}) en ne conservant que les objets dont l'indice est~$\leq i$.

\begin{lemme3} Le $\Z$-module $S_i$ est libre et les éléments $\varepsilon_k\otimes \varepsilon_l$ et $\varepsilon_l \otimes \varepsilon_k\in Y_i\otimes X_i$ ont même image dans $S_i$ pour $k$ et $l$ dans $[e_i+1,r-m_i]$.
\end{lemme3}

\begin{demo} On procède par récurrence sur $i$. Le cas $i=0$ résulte du lemme~\ref{lemmeEtCran}. Supposons les assertions vérifiées pour $S_{i-1}$. Le $\Z$-module $S_i$ est la limite inductive du diagramme
$$ S_{i-1} \longleftarrow Y_{i-1}\otimes X_i \longrightarrow S'_i\: .$$
La flèche $Y_{i-1}\otimes X_{i} \rightarrow Y_{i-1}\otimes X_{i-1}$ envoie $\varepsilon_k\otimes \varepsilon_l$ sur $p\: \varepsilon_k\otimes \varepsilon_l$ si $l\in [r-m_i+1,r-m_{i-1}]$, et sur $\varepsilon_k\otimes \varepsilon_l$ sinon. La flèche $Y_{i-1}\otimes X_{i} \rightarrow Y_i\otimes X_i$ envoie $\varepsilon_k\otimes \varepsilon_l$ sur $p \:\varepsilon_k\otimes \varepsilon_l$ si $k\in [e_{i-1}+1,e_i]$, et sur $\varepsilon_k\otimes \varepsilon_l$ sinon. On en conclut  que $S_i$ est libre, et que $\varepsilon_k\otimes \varepsilon_l$ et $\varepsilon_l \otimes \varepsilon_k\in Y_i\otimes X_i$ ont même image dans $S_i$ pour $k$ et $l$ dans $[e_i+1,r-m_i]$, ce qui achève la récurrence.
\end{demo}

\begin{corollaire3}  \label{corollaireLibre} Le $\Z$-module $S_0^w$ est libre et le groupe $E_0^w$ est un tore.
\end{corollaire3}

\subsubsection{Oubli du niveau}

Le morphisme d'oubli du niveau $\mathcal{M}^w_0\rightarrow \mathcal{M}$ se factorise en un morphisme $\mathcal{B}_0^w \rightarrow \mathcal{B}$ et en un morphisme entre torseurs équivariant sous l'isogénie $E_0^w\rightarrow E$. Les flèches $\mathcal{M}_0^w \rightarrow \mathcal{M}$ et $\mathcal{M}_0\rightarrow \mathcal{M}$ sont {propres}.
Soit $M$ le $1$-motif universel sur $\mathcal{M}$. Notons $M[p]$ le groupe fini et plat sur $\mathcal{M}$ qui est associé à $M$. L'espace de modules relatif $\mathcal{M}_0\rightarrow \mathcal{M}$ paramètre alors les drapeaux isotropes de type $\mathcal{D}$ de $M[p]$.

\subsubsection{Variante de l'espace de modules avec structure parahorique} \label{partieVariante}
On fixe $w\in \mathcal{W}$ et un diagramme commutatif de $\Z$-modules 
\begin{equation}\label{equationXYgras}
\xymatrix{ \mathbb{Y}_0 \ar[r] \ar[d] & \mathbb{Y}_1 \ar[r] \ar[d] & \cdots  \ar[r] & \mathbb{Y}_s \ar[d] \\
\mathbb{X}_0 \ar[d]  & \mathbb{X}_1 \ar[l]\ar[d] & \cdots \ar[l] & \mathbb{X}_s\ar[d] \ar[l] \\
\mathbb{Y}_0 \ar[r] & \mathbb{Y}_1 \ar[r] & \cdots \ar[r]  & \mathbb{Y}_s }
\end{equation}
qui vérifie les conditions de la proposition~\ref{propParaSym}. On note
$$\mathcal{M}^w_{\mathbb{X}_\bullet,\mathbb{Y}_\bullet,\: 0}$$
le champ de modules des $1$-motifs de genre $g$, principalement polarisés, munis d'une structure de niveau parahorique de type $\mathcal{D}$ et d'invariant $w$ en $p$, et munis d'un isomorphisme du diagramme~(\ref{equationXY}) (obtenu grâce à la proposition~\ref{propParaSym}) sur le diagramme~(\ref{equationXYgras}). Tout isomorphisme entre
$$\xymatrix{\mathbb{Y}_0\ar[r] &  \mathbb{X}_0 \ar[r] & \mathbb{Y}_0 } \:\: \mathrm{et} \:\: \xymatrix{\mathbb{Y} \ar[r]^{\times \: p} & \mathbb{X} \ar[r]^{\mathrm{Id}} &\mathbb{Y}}\: ,$$
induit une immersion ouverte et fermée 
$$\mathcal{M}^w_{0,\: \mathbb{X}_\bullet,\mathbb{Y}_\bullet}\longrightarrow \mathcal{M}^w_0\: .$$
Les champs de modules avec $\mathbb{X}_\bullet$ et $\mathbb{Y}_\bullet$ fixés serviront lors de la description précise de la stratification des compactifications (théorème~\ref{thPrinc}).

\section{Cartes locales}

\subsection{Compactifications toroïdales partielles}

On fixe un $\Z$-module libre $X$ de rang $g$ et l'on note $B(X)=\mathrm{Hom}({\mathrm{Sym}}^2(X),\Z)$
le $\Z$-module des applications bilinéaires symétriques entières sur $X$, et
$C(X) \subset B(X)\otimes \R$
le cône des applications bilinéaires symétriques semi-définies positives à radical rationnel. On se donne une décomposition polyédrale $\mathrm{GL}(X)$-admissible $\Sigma$ de $C(X)$ qui est lisse pour la structure entière $B(X)$ (paragraphe IV.2 de~\cite{Deg@FaltingsChai}).
Soit $< \: , \: >$ désigne l'accouplement de dualité entre $B(X\otimes \R)$ et $\mathrm{Sym}^2(X\otimes \R)$. On pose
\begin{eqnarray*}
\sigma^{\vee} & = & \lbrace \: x \in \mathrm{Sym}^2(X\otimes \R) \: \mathrm{tq} \: <\sigma,x>\: \geq \:0 \: \rbrace \\
\sigma^{\perp} & = & \lbrace \: x \in \mathrm{Sym}^2(X\otimes \R) \: \mathrm{tq} \: <\sigma,x>\: = \: 0 \: \rbrace \
\end{eqnarray*}
pour tout $\sigma\in \Sigma$.
Ce sont respectivement un monoïde et un groupe pour l'addition.
Pour tout $\sigma\in \Sigma$, on désigne par $X_{\sigma}$ le plus petit quotient de $X$ tel que toutes les formes bilinéaires $b\in\sigma$ se factorisent par $X_{\sigma} \otimes X_{\sigma}$, et l'on note $r$ son rang. Par exemple si
$\sigma$ est dans l'intérieur de $C(X)$, on a $X_{\sigma}=X$ et si $\sigma=\{0\}$, on a $X_{\sigma}=0$.
On pose
$$S_\sigma= \mathrm{Sym}^2(X_{\sigma}) \: \: \mathrm{et} \:  \: \overline{E}_{\sigma} = \Spec(\Z[S_{\sigma}\cap \sigma^{\vee}]) \: .$$
Ce dernier schéma est un {plongement torique} affine,  muni d'une stratification indexée par les faces de $\sigma$~:~la face $\tau$ correspond à la strate $E_{\sigma}^{\tau}=\Spec(\Z[S_{\sigma} \cap \tau^{\perp}])$ qui est un tore quotient de $E_{\sigma}$. Par exemple, l'origine $\lbrace 0 \rbrace$ correspond à l'unique strate ouverte $E_{\sigma}$ et le cône $\sigma$ à l'unique strate fermée.  La paramétrisation de la stratification est décroissante, c'est-à-dire que si $\tau$ est une face de $\tau'$ la strate attachée à $\tau'$ est  incluse dans l'adhérence de la strate attachée à $\tau$.

On définit un faisceau constructible $\underline{X}_{\sigma}$ sur $\overline{E}_\sigma$, qui est un quotient du faisceau constant $X_{\overline{E}_\sigma}$, en donnant sa restriction aux strates~:
$$\underline{X}_{\sigma}|_{E_{\sigma}^{\tau}}=X_{\tau}\: .$$
La restriction de $\underline{X}_{\sigma}$ à la strate ouverte est nulle et sa restriction à la  strate $E_{\sigma}^{\tau}$ est égale à $X_{\sigma}$ si $\tau$ est inclus dans l'intérieur de $C(X_{\sigma})$.
Par construction, $\mathrm{Sym}^2(X_{\sigma})$ se plonge dans l'anneau des fonctions régulières sur $E_{\sigma}$, donc dans le corps des fractions rationnelles de $\overline{E}_\sigma$. En considérant le diviseur associé à une fraction rationnelle, on obtient une forme bilinéaire symétrique ${B}_{\sigma} \: : \: {X}_{\sigma} \otimes {X}_{\sigma} \rightarrow \underline{\mathrm{Div}}_{\overline{E}_\sigma}$ à valeurs dans le faisceau des diviseurs de $\overline{E}_\sigma$.

\begin{proposition2} \label{propBXxi} La forme $B_\sigma$ se factorise en une forme bilinéaire symétrique
$$\underline{B}_{\sigma} \: : \: \underline{X}_{\sigma} \otimes \underline{X}_{\sigma} \longrightarrow \underline{\mathrm{Div}}_{\overline{E}_\sigma}$$
telle que pour toute section locale $x$ de $\underline{X}_{\sigma}$ le diviseur $\underline{B}(x,x)$ soit effectif à support égal au lieu de non annulation de $x$.
\end{proposition2}

\begin{demo} Vérifions la factorisation sur chaque strate $E_{\sigma}^{\tau}$, où $\tau$ est une face de $\sigma$. Notons $X^{\tau}$ le noyau de la surjection $X \rightarrow X_{\tau}$. Il suffit de montrer que la fonction associée à tout élément $x\otimes y \in X^{\tau}\otimes X_\sigma$ est sans zéro ni pôle sur $E_{\sigma}^{\tau}$. Par construction, on a
$$\mathrm{H}^0(E_{\sigma}^{\tau},\mathcal{O}_{E_{\sigma}^{\tau}}) = \Z[\: s\in X_{\sigma}\otimes X_{\sigma} \:\: \mathrm{tq}  \:\: b(s)=0 \: \: \forall b \in \tau] \: . $$
Or si $b\in \tau$ on a  $b(x\otimes y)=0$, donc $x\otimes y$ induit une fonction régulière sur $E_{\sigma}^{\tau}$. C'est également le cas de la fonction inverse qui est associée à $x \otimes (-y)$, d'où le résultat.

Montrons la seconde assertion. Le cône dual de $C(X_\sigma)$ est $C(X_\sigma^\vee)$ d'où une inclusion $\sigma^\vee \supset C(X_\sigma^\vee)$. Si $x$ est une section locale de $X_\sigma$, on voit $x\otimes x$ comme un élément de $C(X_\sigma^\vee)$. La fonction associée à $x \otimes x$ est régulière sur $\overline{E}_\sigma$ et d'après les principes généraux sur les variétés toriques, le lieu de ses zéros est réunion des strates associées aux faces $\tau$ sur lesquelles $x\otimes x$ prend des valeurs strictement positives comme forme linéaire sur $B(X_\sigma)\otimes \R$.
\end{demo}

On appelle $\Xsi$ le champ de modules des $1$-motifs de genre $g$ principalement polarisés et rigidifiés par $X_{\sigma}$. Il ne dépend que de $X_{\sigma}$ et pas de $\sigma$. On pose également $\mathcal{B_{\sigma}}=\Hom(X_{\sigma},\Agr)$, où $r$ est le rang de $X_\sigma$. Dans le chapitre précédent, on a montré que $\Xsi$ est un $E_{\sigma}$-torseur sur $\mathcal{B}_{\sigma}$. Notons $\Xsib$ le produit contracté
$$\Xsib = \Xsi \times^{E_{\sigma}} \overline{E}_{\sigma} \: .$$

Le champ algébrique $\Xsib$ est affine sur $\mathcal{B}_{\sigma}$ et  muni d'une stratification  indexée par les faces de $\sigma$. On note $Z_{\sigma}^{\tau}$ la strate associée à la face $\tau$ et $Z_{\sigma}=Z_{\sigma}^{\sigma}$ l'unique strate fermée.  Le faisceau $\underline{X}_\sigma$ et la forme bilinéaire $\underline{B}_\sigma$ s'étendent à $\Xsi \times_{\mathcal{B}_{\sigma}} \overline{E}_\sigma$ de manière $E_\sigma$-invariante et se descendent à $\Xsib$. On note toujours $\underline{X}_\sigma$ le faisceau constructible  et 
$$\underline{B}_\sigma \: : \: \underline{X}_\sigma \otimes \underline{X}_\sigma \longrightarrow \underline{\mathrm{Div}}_{\Xsib}$$
la forme bilinéaire obtenus sur $\Xsib$. Le diviseur $\underline{B}_\sigma(x\otimes x)$ est effectif de support le lieu de non annulation de $x\in \underline{X}_\sigma$.
Si $X_\sigma=X$ alors on a $r=g$, on a $\Xsi=E_\sigma$ et l'on a $\Xsib=\overline{E}_\sigma$. Si $\sigma=\lbrace 0 \rbrace$, on a $r=0$ et $\Xsib=\Xsi=\ag$.

\smallskip

\begin{remarque2} Nous ne considérerons des plongements toriques non affines que dans l'énoncé du théorème principal, pour obtenir une description globale du voisinage de l'infini des compactifications. Dans la construction des compactifications et plus précisément dans les propositions~\ref{propRbEt} et~\ref{propR0bet} nous utiliserons le fait que les cônes de $\Sigma$ ne se recouvrent pas.
\end{remarque2}

\subsubsection{Interprétation modulaire à la frontière}

Le champ $\Xsib$ n'admet pas d'interprétation modulaire explicite. En revanche, l'immersion ouverte $(\Xsi\hookrightarrow \Xsib)$ en admet une en tant que \og catégorie fibrée à la frontière \fg dans le sens suivant.
Soient $F$ et $G$ deux champs sur $\Spec(\Z)$ et $F\hookrightarrow G$ une immersion ouverte. 
\begin{definition3} La catégorie fibrée à la frontière $(F\hookrightarrow G)$ est la catégorie fibrée en groupoïdes sur $\mathsf{Ouv}$ qui associe à tout objet  $(U\hookrightarrow S)$  la catégorie des diagrammes $2$-commutatifs 
$$\xymatrix{U \ar[r] \ar[d] & F \ar[d] \\
S \ar[r] & G}$$
\end{definition3}
On appelle $(U\hookrightarrow S)$\textit{-point à la frontière} de $(F\hookrightarrow G)$ la catégorie des diagrammes précédents~;~on la note $(F\hookrightarrow G)(U\hookrightarrow S)$.
 
\smallskip

\begin{remarque3} La catégorie $\mathsf{Ouv}$ forme un site si on la munit de la topologie induite par la topologie \textit{fpqc} sur $U$ et $S$, et la catégorie fibrée $(F\hookrightarrow G)$ est un champ sur ce site.
\end{remarque3}
 
Donnons à présent une description explicite de la catégorie fibrée à la frontière $$(\Xsi\hookrightarrow\Xsib)\: .$$
Soient $(U\hookrightarrow S)$ un objet de $\mathsf{Ouv}$, et $D$ un diviseur effectif irréductible réduit de $S$. On note $p_D\: : \: \underline{\mathrm{Div}}_S \rightarrow \Z_D$ le morphisme qui associe à tout diviseur sa multiplicité le long de $D$. On rappelle que l'on a associé un faisceau constructible $\underline{X}$ et une forme bilinéaire symétrique $\underline{B}:\underline{X}\otimes \underline{X}\rightarrow \underline{\mathrm{Div}}_S$ à tout $1$-motif de Mumford principalement polarisé sur $(U\hookrightarrow S)$.

\begin{proposition3} \label{propmodfrontXsib} $(\Xsi\hookrightarrow \Xsib)(U\hookrightarrow S)$ est le groupoïde des $1$-motifs de Mumford sur $(U\hookrightarrow S)$ de genre $g$, principalement polarisés, rigidifiés  par $X_{\sigma}$ sur $U$, tels que pour toute section locale $x$ de $\underline{X}$, le diviseur $\underline{B}(x,x)$ soit effectif et pour tout diviseur effectif irréductible réduit $D$ de $S$, l'application bilinéaire symétrique $p_D \circ \underline{B} \: : \: X_{\sigma} \otimes X_{\sigma} \rightarrow \Z$ soit incluse dans l'adhérence de $\sigma$.
\end{proposition3}

\begin{demo} Il est équivalent de se donner la flèche composée $$S\longrightarrow \Xsib \longrightarrow \mathcal{B}_\sigma$$ ou de se donner un schéma semi-abélien $\tilde{G}$ de rang torique constant sur $S$, dont le groupe des caractères de la partie torique est $X_\sigma$ et dont la partie abélienne est principalement polarisée de genre $g-r$. Il est  équivalent de se donner un relèvement $U\rightarrow \Xsi$ et un morphisme $X_\sigma\rightarrow \tilde{G}$ sur $U$ induisant un $1$-motif de Mumford principalement polarisé sur $(U\hookrightarrow S)$.

Il suffit de vérifier que sous les conditions données sur $\underline{B}$, la flèche $U\rightarrow \Xsi$ s'étend en une flèche $S\rightarrow \Xsib$. Une éventuelle extension est unique par densité de $U$ dans $S$. Par descente pour la topologie de Zariski, on peut donc supposer que $S=\Spec(R)$ est un schéma local intègre de point générique $\Spec(K)$ et que le torseur $\Xsi \rightarrow \mathcal{B}_\sigma$ est trivial. On obtient ainsi une flèche $U\rightarrow E_\sigma$, que l'on veut étendre en une flèche $S\rightarrow \overline{E}_\sigma$.
On veut montrer que la flèche $\mathrm{H}^0({E}_{\sigma},\mathcal{O})=\Z[\mathrm{Sym}^2(X_{\sigma})]\rightarrow K$  envoie $\mathrm{H}^0(\overline{E}_{\sigma},\mathcal{O})$ dans $R$. Par normalité de $S$, cela se teste en les idéaux premiers de hauteur $1$. Il suffit de voir que pour tout $x\otimes y \in \mathrm{H}^0(\overline{E}_{\sigma},\mathcal{O}) \subset \Z[\mathrm{Sym}^2(X_{\sigma})]$ et tout diviseur $D$ de $S$, on a $p_D \circ \underline{B}(x\otimes y) \geq 0$. C'est évident car 
$$\mathrm{H}^0(\overline{E}_{\sigma},\mathcal{O})=\Z[\: \lbrace u\otimes v\in X_{\sigma}\otimes X_{\sigma} \: \: \mathrm{tq} \: \: \forall b \in \sigma, \: b(u,v)\geq 0 \rbrace\: ]$$
et la forme bilinéaire $p_D \circ \underline{B}$ est incluse dans l'adhérence de $\sigma$ par hypothèse.
\end{demo}

\begin{remarque3} Si $S$ est local, la flèche $S\rightarrow \Xsib$ envoie le point fermé de $S$ dans la strate $Z_{\sigma}$ si et seulement si $\sigma$ est  minimal parmi les cônes de $\Sigma$ dont l'adhérence contient la forme bilinéaire $p_D \circ \underline{B}$ pour tout diviseur $D$ de $S$.
\end{remarque3}

\subsubsection{Cas d'une structure de niveau parahorique} \label{partieCompactPara}

On considère encore la décomposition $\Sigma$ du paragraphe précédent, qui est lisse pour la structure entière $B(X)$. Soit $\sigma\in \Sigma$. On note $X_{\sigma}$ le quotient de $X$ associé  à $\sigma$, et $r$ son rang.
On note $\XsiK$ le champ algébrique des $1$-motifs de genre $g$ qui sont principalement polarisés, rigidifiés par $X_{\sigma}$, et munis d'une structure de niveau parahorique de type $\mathcal{D}$ en $p$. On a vu dans le chapitre précédent que $\XsiK$ admettait une décomposition en parties ouvertes et fermées
$$\XsiK = \coprod_{w\in W_{\sigma}} \Xsiw$$
indexée par le quotient de groupes de Weyl 
$$\mathcal{W}_\sigma =  W_{L} \: \backslash\ W_{\gsp} \: / \ W_{L^\mathcal{D}} \: .$$
Pour tout $w\in W_{\sigma}$, on a décrit dans la proposition \ref{propXi0w} un tore $E_{\sigma,\: 0}^w$ et un champ $\mathcal{B}_{\sigma,\: 0}^w$ tels que $\Xsiw \rightarrow \mathcal{B}_{\sigma,\: 0}^w$ soit un torseur sous $E_{\sigma,\: 0}^w$. Le $\Z$-module libre  $S_{\sigma,\: 0}^w = \Hom(E_{\sigma,\: 0}^w,\Gm)$ définit un plongement torique affine  $E_{\sigma,\: 0}^w \hookrightarrow \overline{E}_{\sigma,\: 0}^w$.
On définit le champ algébrique $\Xsiwb \rightarrow \mathcal{B}_{\sigma,\: 0}^w$ comme le produit contracté
$$\Xsiwb = \Xsiw \times^{E_{\sigma,\:0}^w} \overline{E}_{\sigma,\:0}^w$$
et le champ algébrique $\XsiKb\rightarrow \Spec(\Z)$ comme l'union disjointe
$$\XsiKb=\coprod_{w\in \mathcal{W}_{\sigma}} \Xsiwb \: .$$
Le champ $\XsiKb$ est muni d'une stratification localement fermée de strates
$$Z_{\sigma,\: 0}^{\tau}= \coprod_{w\in \mathcal{W}_{\sigma}} Z_{\sigma,\: 0}^{\tau, \: w}$$
indexées par les faces $\tau$ de $\sigma$.  On note  $Z_{\sigma,\: 0}=Z_{\sigma,\: 0}^{\sigma}$ l'unique strate fermée.

On a une interprétation de $(\XsiKb\hookrightarrow \XsiK)$ comme catégorie fibrée à la frontière analogue à celle de $(\Xsi\hookrightarrow \Xsib)$. Soit $(U\hookrightarrow S)$ un objet de $\mathsf{Ouv}$. Tout $1$-motif de Mumford $M$ sur $(U\hookrightarrow S)$ muni d'une structure de niveau parahorique $H_\bullet$ d'invariant $w\in \mathcal{W}$  définit une famille compatible de formes bilinéaires associées aux $1$-motifs de Mumford $M/H_i$~, donc un élément $\underline{B}(M,H_\bullet)$ de $\mathrm{Hom}(S_{\sigma,\: 0}^w\: , \underline{\mathrm{Div}}_S)$.

\begin{proposition3} \label{propmodfrontXsiKb}  $(\XsiKb\hookrightarrow \XsiK)(U\hookrightarrow S)$ est le groupoïde des $1$-motifs de Mumford sur $(U\hookrightarrow S)$ de genre $g$, principalement polarisés, rigidifiés par $X_{\sigma}$ et munis d'une structure de niveau parahorique de type $\mathcal{D}$ en $p$, tels que pour toute section locale $x$ de $X_\sigma$ le diviseur $\underline{B}(M,H_\bullet)(x,x)$ soit effectif et pour tout diviseur effectif irréductible réduit $D$ de $S$, la forme bilinéaire symétrique $p_D \circ \underline{B}(M,H_\bullet)$ soit incluse dans l'adhérence de $\sigma$.
\end{proposition3}

\subsubsection{Oubli du niveau}

L'inclusion des structures entières  $S_\sigma\hookrightarrow S_{\sigma,\: 0}^w$ donne une flèche propre $$\overline{E}_{\sigma,\: 0}^w \rightarrow \overline{E}_{\sigma}$$ et la fonctorialité du produit contracté donne  une flèche propre $$\Xsiwb \longrightarrow \Xsib$$ qui étend l'oubli du niveau $\Xsiw \rightarrow \Xsi$. Ces morphismes induisent un diagramme cartésien 
$$
\xymatrix{ \XsiK \ar[r] \ar[d] & \XsiKb \ar[d] \\
\Xsi \ar[r] & \Xsib}
$$
dont les flèches verticales sont propres et dont les flèches horizontales sont des unions disjointes de plongements toroïdaux.

\begin{remarque3} Comme la décomposition $\Sigma$ est lisse pour la structure entière $B(X)$, le champ $\Xsib$ est lisse sur $\Spec(\Z)$ et $\Xsi\hookrightarrow \Xsib$ est le complémentaire d'un diviseur à croisements normaux. Le produit fibré $$\XsiK \times_{\Spec(\Z)} \Spec\left(\Z[1/p]\right)$$ est lisse sur $\Spec(\Z[1/p])$ car les champs $\agr$ sont lisses sur $\Spec(\Z[1/p])$. Toutefois, comme $\Sigma$ n'est pas lisse pour les différentes structures entières $B_{\sigma,\: 0}^w$, le champ $\Xsiwb$ n'est pas lisse sur $\Spec(\Z[1/p])$ et $$\Xsiw \hookrightarrow \Xsiwb$$ n'est pas le complémentaire d'un diviseur à croisements normaux. De même, les singularités de $$\Xsiwb\longrightarrow \Spec(\Z_p)$$ ne sont pas uniquement celles des $\agrK$ avec $r\leq g$.
\end{remarque3}

\subsection{Cartes locales formelles}

On choisit une présentation étale surjective $\mathcal{U} \rightarrow \mathcal{B}_{\sigma}$
du champ algébrique $\mathcal{B}_{\sigma}$ par un schéma affine ne dépendant que du rang de  $X_\sigma$. On note 
$${\widehat{\overline{\mathcal{M}}}_{\sigma}} = (\widehat{\mathcal{U}\times_{\mathcal{B}_{\sigma}} \Xsib})^{\mathcal{U}\times_{\mathcal{B}_{\sigma}} Z_\sigma}$$
le complété du schéma affine $\mathcal{U}\times \Xsib$ le long du sous-schéma fermé $\mathcal{U}\times Z_\sigma$. Ce schéma dépend du choix de $\mathcal{U}$, même si on omet ce dernier des notations pour simplifier. Il contient le sous-schéma ouvert
$${\widehat{\mathcal{M}}_{\sigma}} =   \Xsi \times_{\Xsib} {\widehat{\overline{\mathcal{M}}}_{\sigma}}\: ,$$
le sous-schéma fermé de complétion $\mathcal{U}\times Z_\sigma$ que l'on note encore $Z_\sigma$ et le diviseur à croisements normaux
$$\partial {\widehat{\overline{\mathcal{M}}}_{\sigma}}={\widehat{\overline{\mathcal{M}}}_{\sigma}}-{\widehat{\mathcal{M}}_{\sigma}} \: .$$
On a une immersion fermée $Z_{\sigma}\hookrightarrow \partial {\widehat{\mathcal{M}}_{\sigma}}$.
On note $M_\sigma$ le $1$-motif principalement polarisé universel sur $\Xsi$. Il induit un $1$-motif de Mumford encore noté $M_{\sigma}$ sur $$\left(\widehat{\Xsi}\hookrightarrow \widehat{\overline{\mathcal{M}}}_{\sigma}\right)\: .$$
Le schéma ${\widehat{\overline{\mathcal{M}}}_{\sigma}}$ est normal, affine, complet relativement à $Z_\sigma$, excellent, et par construction, $M_{\sigma}$ est positif relativement à $Z_{\sigma}$. Les hypothèses de la construction de Mumford sont vérifiées donc il existe un schéma semi-abélien
$$G_{\sigma}\longrightarrow {\widehat{\overline{\mathcal{M}}}_{\sigma}}$$
qui est un objet de la catégorie $$\mathsf{Deg}_{\widehat{\Xsi},\widehat{\overline{\mathcal{M}}}_{\sigma},Z_\sigma,\: pol}\: .$$
On lui associe un faisceau constructible $\underline{X}(G_\sigma)$ par le théorème~\ref{thFaisCons}  et une application bilinéaire
$$\underline{B}(G_\sigma) \: : \: \underline{X}(G_\sigma)\otimes \underline{X}(G_\sigma)\longrightarrow \mathrm{Div}_{\widehat{\overline{\mathcal{M}}}_{\sigma}}$$
par la proposition~\ref{propBGsab}. On rappelle que l'on a également construit un faisceau constructible $\underline{X_\sigma}$ et une forme bilinéaire $\underline{B}_\sigma$ sur $\Xsib$. On note de la même manière les objets obtenus par image inverse sur $\widehat{\overline{\mathcal{M}}}_{\sigma}\: .$

On appelle $\tilde{G}_{\sigma}$ le schéma semi-abélien sur $\widehat{\overline{\mathcal{M}}}_{\sigma}$ qui provient du schéma semi-abélien de rang torique constant universel sur $\mathcal{B}_\sigma$. Il est égal à la partie semi-abélienne de $M_\sigma$.
Rappelons le résultat suivant, dû à Faltings et Chai.

\begin{proposition2}[\cite{Deg@FaltingsChai}, proposition IV.3.1] \label{propMumUniv} Il existe des isomorphismes canoniques $\underline{X}(G_\sigma)\iso \underline{X}_\sigma$ et $\underline{B}(G_\sigma)\iso \underline{B}_\sigma$. La restriction de $G_\sigma$ à $\widehat{\Xsi}$ est un schéma abélien principalement polarisé de genre $g$. La restriction de $G_\sigma$ à $Z_\sigma$ est canoniquement isomorphe à $\tilde{G}_\sigma$.  L'application de Kodaira-Spencer réalise un isomorphisme
$$K_{G_{\sigma}/ \: {\widehat{\overline{\mathcal{M}}}_{\sigma}}} \: : \: \mathrm{Sym}^2 (\Omega_{G_{\sigma}}) \: \isolong \: \Omega^1_{{\widehat{\overline{\mathcal{M}}}_{\sigma}}/\Z}\left(\mathrm{log} \:  \partial {\widehat{\overline{\mathcal{M}}}_{\sigma}}\right) \: .$$
\end{proposition2}

Le point délicat est de montrer la bijectivité de la flèche de Kodaira-Spencer obtenue par la proposition \ref{propKSMum}. Cette flèche est construit \textit{via} le $1$-motif $M_\sigma$ et la bijectivité provient de l'universalité de $M_\sigma$.

\begin{remarque2} Il est crucial d'utiliser la {lissité} de $\Xsi \rightarrow \Spec(\Z)$ pour montrer que la flèche de Kodaira-Spencer réalise un isomorphisme. Cet isomorphisme nous sera indispensable pour contrôler le processus d'approximation.
\end{remarque2}

On  dispose par conséquent d'une flèche ${\widehat{\mathcal{M}}_{\sigma}}\rightarrow \ag$. On  dit que le diagramme
$$\left(\xymatrix{{\widehat{\overline{\mathcal{M}}}_{\sigma}} & {\widehat{\mathcal{M}}_{\sigma}} \ar[l] \ar[r] & \ag}\right)$$
est la \textit{carte locale formelle sans niveau} associée à $\sigma\in\Sigma$.

\subsubsection{Propriété modulaire proche de la frontière}

Nous développons un formalisme analogue à celui des espaces de modules à la frontière en tenant compte de la nature complète des anneaux intervenant dans la construction de Mumford. Nous cherchons  notamment à imposer des conditions relatives à un sous-champ algébrique fermé.

Pour tout schéma $T$, on définit la catégorie $\mathsf{OuvComp}_T$ comme la sous-catégorie de $\mathsf{OuvComp}$ dont les objets se fibrent sur $T$~;~les flèches sont celles qui induisent l'identité sur $T$. On se donne un diagramme de champs sur $T$
$$F\longrightarrow G \longleftarrow Z$$ 
tel que $F\rightarrow G$ soit une immersion ouverte, $Z\rightarrow G$ soit une immersion fermée et $F\times_G Z$ soit vide, c'est-à-dire que $F$ et $Z$ soient disjoints dans $G$.

\begin{definition3} La catégorie fibrée proche de la frontière $(F\rightarrow G\leftarrow Z)_T$ est la catégorie fibrée en groupoïdes sur $\mathsf{OuvComp}_T$ qui associe à tout objet  $(U\rightarrow S \leftarrow S_0)$ la catégorie $$(F\longrightarrow G \longleftarrow Z)_T (U\longrightarrow S \longleftarrow S_0)$$ des diagrammes $2$-commutatifs 
$$\xymatrix{U \ar[r] \ar[d] & F \ar[d] \\
S \ar[r] & G\\
S_0 \ar[r] \ar[u] & Z \ar[u] }$$
compatibles aux projections $G\rightarrow T$ et $S\rightarrow T$.
\end{definition3}

Soit $(U\rightarrow S \leftarrow S_0)$ un objet de $\mathsf{OuvComp}_T$. On appelle $(U\rightarrow S \leftarrow S_0)$-point proche de la frontière de $(F\rightarrow G \leftarrow Z)_T$ la catégorie des diagrammes précédents.
On donne une description explicite de la catégorie fibrée proche de la frontière
$$(\: \widehat{\Xsi}\: \longrightarrow\: {\widehat{\overline{\mathcal{M}}}_{\sigma}} \:\longleftarrow\: Z_\sigma\: )_\mathcal{U} \: .$$
On note comme précédemment $\tilde{G}_\sigma$ le schéma semi-abélien de rang torique constant universel sur $\mathcal{B}_\sigma$. Il induit un schéma semi-abélien $\tilde{G}_{\sigma}$ de rang torique constant sur $S_0$ par image inverse~\textit{via}
$$S_0 \longrightarrow \mathcal{U} \longrightarrow \mathcal{B}_\sigma \: .$$
Si $G$ est un schéma semi-abélien sur $S$, tout isomorphisme $G_{S_0}\iso \tilde{G}_\sigma$ sur $S_0$ induit une surjection $X_\sigma \rightarrow \underline{X}(G)$ de faisceaux constructibles sur $S$.

\begin{proposition3} \label{propModProchFront} Le foncteur qui associe à $f:S\rightarrow \widehat{\overline{\mathcal{M}}}_{\sigma}$ le schéma semi-abélien $f^*G_\sigma$ réalise une équivalence de catégories entre
$$(\widehat{\mathcal{M}}_{\sigma} \longrightarrow \widehat{\overline{\mathcal{M}}}_{\sigma} \longleftarrow Z_{\sigma})_\mathcal{U}(U\longrightarrow S \longleftarrow S_0)$$
et la catégorie des couples $(G,\:G|_{S_0} \iso \tilde{G}_\sigma)$, où $G$ est un objet de $\mathsf{Deg}_{U,S,S_0,\: pol}$ de dimension relative $g$ principalement polarisé sur $U$ tel que pour tout diviseur effectif irréductible réduit $D$ de $S$, l'application bilinéaire  $p_D \circ \underline{B}(G)  :  X_{\sigma} \otimes X_{\sigma} \rightarrow \Z$ soit incluse dans l'adhérence de $\sigma$, et tel que $\sigma$ soit minimal parmi les cônes de $\Sigma$ vérifiant cette propriété.
\end{proposition3}

\begin{demo} Les schémas semi-abéliens obtenus par image inverse de $G_\sigma$ vérifient les propriétés requises. Inversement, soit $G$ un schéma semi-abélien sur $S$ vérifiant les hypothèses de la proposition. Le  théorème~\ref{ThConsMum}~(construction de Mumford inverse) permet de lui associer un $1$-motif de Mumford sur $(U\hookrightarrow S)$ et la proposition~\ref{propmodfrontXsib} donne un $(U\hookrightarrow S)$-point de $\Xsib$. L'hypothèse de minimalité de $\sigma$ implique que cette flèche envoie $S_0$ dans $Z_{\sigma}$. Comme le schéma affine $S$ est complet relativement à $S_0$, la flèche $S\rightarrow \Xsib$ se factorise en une flèche 
$$S\longrightarrow \widehat{\overline{\mathcal{M}}}_{\sigma} \: .$$
La fonctorialité de la construction de Mumford montre que $G$ est l'image inverse de $G_\sigma$ par cette flèche.
\end{demo}

\subsubsection{Algébricité du groupe $p$-torsion}

Le schéma semi-abélien $G_\sigma$ ne  descend pas à $\Xsib$. En revanche, le groupe des points de $p$-torsion $G_\sigma[p]\rightarrow \widehat{\Xsi}$ descend à $\Xsi$. En effet, d'après le théorème~\ref{ThConsMum}, il existe un isomorphisme canonique
$$G_\sigma[p] \isolong M_\sigma[p]$$
sur $\widehat{\Xsi}$. Par construction, le $1$-motif $M_\sigma$ se descend à $\Xsi$ et son groupe de $p$-torsion $M_\sigma[p]$ également.
On résume ceci en disant que $G_\sigma$ est une donnée de nature transcendante, alors que  $G_\sigma[p]$ est algébrique. Cette remarque sera {fondamentale} pour construire une compactification toroïdale de $\agK$.

\subsection{Approximation des cartes locales formelles}

On vient de dire que le schéma semi-abélien 
$$G_\sigma \longrightarrow \widehat{\overline{\mathcal{M}}}_{\sigma}$$
ne descend pas à $\Xsib$. En employant des techniques d'Artin, on peut cependant l'{approcher} par un schéma  $G^{h}_\sigma$ qui est semi-abélien sur un schéma étale sur $\Xsib$. La nouveauté la plus significative par rapport à~\cite{Deg@FaltingsChai} est que l'on effectuera cette approximation en conservant un isomorphisme $G_\sigma^h[p] \iso M_\sigma[p]$, ce qui sera indispensable pour notre construction des compactifications avec niveau.

Une différence mineure par rapport à~\cite{Deg@FaltingsChai} réside dans l'emploi de la notion de couple hensélien et du théorème de Popescu. Cela évite de localiser les cartes locales formelles en tous leurs points fermés~(\cite{Deg@FaltingsChai} après la définition~IV.4.5). 

\subsubsection{Couples henséliens} Rappelons rapidement cette notion.

\begin{definition3}[\cite{Hens@Gabber}] Soit $A$ un anneau et $I$ un idéal. On dit que le couple $(\Spec(A)\: ,\: \Spec(A/I))$ est hensélien si une des conditions équivalentes suivantes est vérifiée~:
\begin{itemize}
\item pour tout $f\in A[t]$ et toute racine simple $\bar{a}$ de $f$ dans $A/I$, il existe une racine $a\in A$ de $f$ relevant $\bar{a}$,
\item pour tout diagramme commutatif de schémas
\begin{eqnarray*}
\xymatrix{S \ar[rr]^{etale} & & \Spec(A) \\
 & \Spec(A/I) \ar[lu]_{p} \ar[ru] & }
\end{eqnarray*}
le morphisme $p$ se relève en un morphisme de $\Spec(A)$-schémas $\Spec(A)\rightarrow S$.
\end{itemize}
\end{definition3}

En particulier, un schéma local $S$ de point fermé $s$ est hensélien si et seulement si le couple $(S,s)$ est hensélien.

\begin{definition3} Soient $A$ un anneau et $I$ un idéal. L'hensélisé de $$\left(\Spec(A)\: , \: \Spec(A/I)\right)$$ est le couple $$\left(\Spec(A^h)\: ,\: \Spec(A^h/I^h)\right)$$ défini par $$\Spec(A^h) = \varprojlim S \: \: \mathrm{et} \: \: I^h=IA^h \: ,$$
où la limite projective est prise sur la catégorie des diagrammes commutatifs précédents. On a $A^h/I^h \simeq A/I$.
\end{definition3}

On peut donc définir l'hensélisé $S^{h_{S'}}$ d'un schéma affine~$S$ le long d'un sous-schéma fermé~$S'$. 

\begin{remarque3} La notion de couple strictement hensélien n'a pas de sens, ni celle d'hensélisé strict le long d'un sous-schéma fermé. On ne peut pas définir l'hensélisé d'un schéma non affine le long d'un sous-schéma fermé car les limites projectives n'existent pas toujours dans la catégorie des schémas.
\end{remarque3}

\subsubsection{Approximation selon Popescu}
 
Un morphisme de schémas affines $\Spec(B)\rightarrow \Spec(A)$ est dit \textit{régulier} s'il est plat à fibres géométriques régulières. Si $A$ est local et excellent, la flèche $\Spec(\widehat{A}) \rightarrow \Spec(A)$ est régulière. Si $A$ est excellent, pour tout idéal $I$ de $A$, le morphisme $$\Spec(\widehat{A}^I) \longrightarrow \Spec(A)$$ du spectre du complété $I$-adique de $A$ dans le spectre de $A$ est régulier par~\egalong{iv}{2}{7.8.3}. Popescu montre le théorème fondamental suivant.

\begin{theoreme3}[\cite{Popescu@Swan}] \label{thPopescu} Une flèche $\Spec(B)\rightarrow \Spec(A)$  entre schémas affines n{\oe}thériens est régulière si et seulement si elle fait de $\Spec(B)$ une limite projective filtrante de schémas affines lisses sur $\Spec(A)$.
\end{theoreme3}
En particulier, si $A$ est un anneau excellent et $I$ un idéal, $\Spec(\widehat{A}^I)$ est limite projective de schémas affines lisses sur $\Spec(A)$.
L'intérêt du théorème de Popescu en ce qui concerne les questions d'approximation résulte du lemme suivant, dont la démonstration a été communiquée par Luc Illusie.

\begin{lemme3} \label{lemApprox} Soient $(A,I)$ un couple hensélien et $n$ un entier naturel. Notons $S=\Spec(A)$ et $S_n=\Spec(A/I^{n+1})$. Soient $S'$ un schéma lisse sur $S$ et $t : S_n \rightarrow S'$ un morphisme de schémas sur $S$. Il existe une section $s:S\rightarrow S'$ telle que les morphismes
$$\xymatrix{S_n \ar[r] &S \ar[r]^s & S'}
\: \: \mathrm{et} \: \: \xymatrix{S_n \ar[r]^t & S'}$$
coïncident.
\end{lemme3}

\begin{demo}
Quitte à remplacer $S'$ par un voisinage ouvert de l'image de $S_n$, on peut supposer qu'il existe une factorisation de la flèche lisse $S' \rightarrow S$ en
\begin{eqnarray*}
\xymatrix{S' \ar[r]^u & \mathbb{A}_S^r \ar[r]^v & S \: ,}
\end{eqnarray*}où $u$ est étale et $v$ est la projection canonique.
On veut trouver une flèche $s:S\rightarrow S'$ qui rende le diagramme suivant commutatif~:

\begin{eqnarray*}
\xymatrix{S' \ar[r]^u & \mathbb{A}_S^r \ar[r]^v & S \ar@{-->}@/_2pc/[ll]_s \\
& S_n \ar[ru] \ar[u] \ar[lu]^t &}
\end{eqnarray*}
Un relèvement $s'  : S \rightarrow \mathbb{A}_S^r$ de $S_n \rightarrow \mathbb{A}_S^r$ existe toujours. Le couple $(S,S_0)$ est hensélien donc~$s'$ se factorise par l'hensélisé de $\mathbb{A}_S^r$ le long de l'image de $$S_0\hookrightarrow S' \longrightarrow \mathbb{A}^r_S\: .$$ Comme $u$ est étale, cet hensélisé est isomorphe à l'hensélisé de $S'$ le long de $S_0$. Cela termine la construction de~$s$.
\end{demo}

Indiquons rapidement comment déduire l'existence d'approximations des résultats précédents. Soient $A$ un anneau n{\oe}thérien excellent, $I$ un idéal et $n$ un entier naturel. Posons $S=\Spec(A)$, $\widehat{S}=\Spec(\widehat{A}^I)$ et $S_n=\Spec(A/I^{n+1})$. Le hensélisé $(S^h,S_0^h)$ de $(S,S_0)$ est excellent en tant que limite projective de schémas de type fini sur $S$. D'après le théorème~\ref{thPopescu}, $\widehat{S}$ est limite projective filtrante de schémas $S^{\alpha}$ lisses sur $S^h$. Toute donnée de type fini sur $\widehat{S}$ provient par image réciproque d'une donnée sur $S^\alpha$ pour $\alpha$ suffisamment grand. D'après le lemme~\ref{lemApprox}, il existe une donnée sur $S^h$ qui a même restriction à $S_n$ que la donnée initiale. On a donc réalisé une approximation à l'ordre~$n$.

\subsubsection{Approximation des cartes locales formelles sans niveau}

Nous allons approcher la carte locale formelle sans niveau associée à $\sigma\in\Sigma$. Définissons le schéma 
$$\overline{\mathcal{M}}_\sigma^h\:=\: {(\mathcal{U}\times_{\mathcal{B_{\sigma}}} \Xsib)}^{h_{\mathcal{U}\times_{\mathcal{B}_{\sigma}} Z_{\sigma}}}\: ,$$
qui est l'hensélisé du schéma affine $\mathcal{U}\times \Xsib$ le long du sous-schéma fermé $\mathcal{U}\times Z_\sigma$. Il est excellent et contient le sous-schéma ouvert 
$$\Xsi^h = \Xsi \times_{\Xsib} \overline{\mathcal{M}}_\sigma^h \:,$$
le diviseur à croisements normaux 
$$\partial \overline{\mathcal{M}}_\sigma^h = \overline{\mathcal{M}}_\sigma^h - \Xsi^h $$
et le sous-schéma fermé d'hensélisation que l'on note encore $Z_\sigma$.

\begin{proposition3} \label{propGxih} Il existe un schéma semi-abélien $G_{\sigma}^h$ sur $\Xsib^h$  vérifiant les propriétés suivantes~:
\begin{itemize}
\item il est abélien principalement polarisé de genre $g$ sur $\Xsi^h$,
\item sa restriction à $Z_\sigma$ est canoniquement isomorphe à $\tilde{G}_\sigma$,
\item $G_{\sigma}^h[p]$ est canoniquement isomorphe à $M_{\sigma}[p]$ sur $\Xsi^h$,
\item $\underline{X}(G_{\sigma}^h)$ est canoniquement isomorphe à $\underline{X}_{\sigma}$,
\item $\underline{B}(G_{\sigma}^h)$ est canoniquement isomorphe à $\underline{B}_{\sigma}$,
\item l'application de Kodaira-Spencer associée à $G_{\sigma}^h$ induit un isomorphisme $$\mathrm{Sym}^2 (\Omega_{G_{\sigma}^h}) \: \isolong  \: \Omega^1_{\Xsib^h/\Z}(\mathrm{log} \: \partial \Xsib^h)\: .$$
\end{itemize}
\end{proposition3}

\begin{demo} D'après le théorème~\ref{thPopescu}, ${\widehat{\overline{\mathcal{M}}}_{\sigma}}$ est limite projective de schémas $\Xsib^{\alpha}$ lisses sur $\Xsib^h$, où $\alpha$ parcourt un ensemble ordonné filtrant. Le schéma semi-abélien $G_\sigma$ provient d'un schéma semi-abélien $G_\sigma^\alpha$ sur $\Xsib^\alpha$ pour $\alpha$ suffisamment grand. Posons
$$\Xsi^\alpha\: = \: \Xsib^\alpha\times_{\Xsib^h} \Xsi^h \quad \mathrm{et} \quad \partial\Xsib^\alpha \: = \: \Xsib^\alpha\times_{\Xsib^h} \partial \Xsib^h$$
Comme $\Xsib^\alpha\rightarrow \Xsib^h$ est lisse, ce sont respectivement un ouvert dense et un diviseur à croisements normaux de $\Xsib^\alpha$. Le schéma $G_\sigma^\alpha$ est abélien sur $\Xsi^\alpha$ et quitte à augmenter $\alpha$, on peut supposer qu'il est principalement polarisé.
On peut également supposer que les isomorphismes 
$$G_\sigma[p]\times {\widehat{{\mathcal{M}}}_{\sigma}} \isolong M_\sigma[p]\: , \: \: G_\sigma\times Z_\sigma \isolong \tilde{G}_\sigma \: , \: \: \underline{X}(G_{\sigma})\isolong \underline{X}_{\sigma} \: \: \mathrm{et} \: \: \underline{B}(G_{\sigma})\isolong \underline{B}_{\sigma}$$
proviennent d'isomorphismes
$$G_\sigma^\alpha[p] \times \Xsi^\alpha \isolong M_\sigma[p]\: , \: \: G_\sigma^\alpha \times Z_\sigma \isolong \tilde{G}_\sigma \: , \: \:  \underline{X}(G_{\sigma}^\alpha)\isolong \underline{X}_{\sigma} \: \: \mathrm{et} \: \: \underline{B}(G_{\sigma}^\alpha)\isolong \underline{B}_{\sigma}\: .$$
Soit $n$ un entier naturel. On applique le lemme~\ref{lemApprox} et l'on obtient une section~$s:\Xsib^h\rightarrow \Xsib^\alpha$ qui relève le morphisme composé
$$Z_{\sigma,\: n} \rightarrow {\widehat{\overline{\mathcal{M}}}_{\sigma}} \rightarrow \Xsib^\alpha \: ,$$
où $Z_{\sigma,\: n}$ est le sous-schéma fermé de $\Xsib^h$ défini par la puissance $n$-ième de l'idéal de $Z_\sigma$.
Il suffit alors de poser $G_\sigma^h = s^* \: G_\sigma^\alpha$ et de vérifier que si $n$ est assez grand, l'application de Kodaira-Spencer associée à $G_\sigma^h$ réalise un isomorphisme, ce qui résulte de~\cite[III.9.4, IV.4.1 et IV.4.2]{Deg@FaltingsChai}.
\end{demo}

Le schéma $G_\sigma^h$ n'est pas unique mais les compactifications ne dépendront pas du choix.

\subsubsection{Recomplétion de l'approximation}

On considère l'image inverse de $G_\sigma^h$ par la flèche 
$${\widehat{\overline{\mathcal{M}}}_{\sigma}}\longrightarrow \Xsib^h$$
et l'on trouve un schéma semi-abélien
$$G_\sigma^h \longrightarrow {\widehat{\overline{\mathcal{M}}}_{\sigma}}\: .$$
D'après la proposition~\ref{propModProchFront}, ce schéma définit une flèche
$$f \: : \: \left( \widehat{\Xsi} \longrightarrow {\widehat{\overline{\mathcal{M}}}_{\sigma}} \longleftarrow Z_\sigma \right) \longrightarrow  \left( \widehat{\Xsi} \longrightarrow {\widehat{\overline{\mathcal{M}}}_{\sigma}} \longleftarrow Z_\sigma \right)$$
de la catégorie $\mathsf{OuvComp}_{\mathcal{U}}$. L'image inverse de $G_\sigma$ par $f$ est isomorphe à $G_\sigma^h$ et $f$ est l'identité au-dessus de $Z_\sigma$. En fait, $f$ est l'identité de 
$$ \left( \widehat{\Xsi} \longrightarrow {\widehat{\overline{\mathcal{M}}}_{\sigma}} \longleftarrow Z_\sigma \right)$$
\textit{modulo} une puissance arbitrairement grande de l'idéal de $Z_\sigma$.

\begin{proposition3} \label{propRecomplAutoEt} Le morphisme $f$ est étale et est un automorphisme  de ${\widehat{\overline{\mathcal{M}}}_{\sigma}}$.
\end{proposition3}

\begin{demo} Les isomorphismes de Kodaira-Spencer associés à $G_\sigma$ et $G_\sigma^h$ définissent un isomorphisme
$$f^* \:  \Omega^1_{{\widehat{\overline{\mathcal{M}}}_{\sigma}}/\Z}(\mathrm{log}\:  \partial {\widehat{\overline{\mathcal{M}}}_{\sigma}})\: \isolong\: \Omega^1_{{\widehat{\overline{\mathcal{M}}}_{\sigma}}/\Z}(\mathrm{log} \: \partial {\widehat{\overline{\mathcal{M}}}_{\sigma}})\: .$$
Comme $f$ est l'identité sur $Z_{\sigma}$, l'isomorphisme précédent induit un isomorphisme
$$f^* \: \Omega^1_{{\widehat{\overline{\mathcal{M}}}_{\sigma}}/\Z} \: \isolong \:  \Omega^1_{{\widehat{\overline{\mathcal{M}}}_{\sigma}}/\Z} \: .$$
Puisque ${\widehat{\overline{\mathcal{M}}}_{\sigma}}\rightarrow \Spec(\Z)$ est formellement lisse, le critère jacobien montre que $f$ est étale. Le schéma $${\widehat{\overline{\mathcal{M}}}_{\sigma}}$$ est hensélien relativement au fermé $Z_{\sigma}$. La restriction de $f$ à ce fermé est un isomorphisme donc $f$ est un automorphisme.
\end{demo}

Nous avons contrôlé l'approximation~:~$f$ est un automorphisme, aucune information n'a été perdue et $G_{\sigma}^h$ est aussi \og universel \fg que $G_{\sigma}$. Soit $(U\rightarrow S \leftarrow S_0)$ un objet de $\mathsf{OuvComp}_{\mathcal{U}}$. La proposition suivante résulte de la proposition~\ref{propModProchFront}.

\begin{proposition3} \label{propModProchFrontHens} Le foncteur qui associe à $f:S\rightarrow \Xsib^h$ le schéma semi-abélien $f^*G^h_\sigma$ réalise une équivalence de catégories entre
$$(\Xsi^h \longrightarrow \Xsib^h \longleftarrow Z_{\sigma})_\mathcal{U}(U\longrightarrow S \longleftarrow S_0)$$
et la catégorie des couples $(G,\:G|_{S_0} \iso \tilde{G}_\sigma)$ où $G$ est un objet de $\mathsf{Deg}_{U,S,S_0,\: pol}$ de dimension relative $g$ qui est principalement polarisé sur $U$, muni  d'un isomorphisme
$G|_{S_0} \iso \tilde{G}_\sigma$
sur $S_0$, tel que pour tout diviseur effectif irréductible réduit $D$ de $S$ l'application bilinéaire $p_D \circ \underline{B} \: : \: X_{\sigma} \otimes X_{\sigma} \rightarrow \Z$ soit incluse dans l'adhérence de $\sigma$, et tel que $\sigma$ soit minimal parmi les cônes de~$\Sigma$ vérifiant cette propriété.
\end{proposition3}

\begin{demo} Comme $S$ est complet relativement à $S_0$, on a
$$(\widehat{\mathcal{M}}_{\sigma} \longrightarrow \widehat{\overline{\mathcal{M}}}_{\sigma} \longleftarrow Z_{\sigma})_\mathcal{U}(U\longrightarrow S \longleftarrow S_0) \: = \: (\Xsi \longrightarrow \Xsib^h \longleftarrow Z_{\sigma})_\mathcal{U}(U\longrightarrow S \longleftarrow S_0) \:.$$
Il suffit  d'appliquer la proposition~\ref{propModProchFront} et de composer la flèche $S\rightarrow \widehat{\overline{\mathcal{M}}}_{\sigma}$
obtenue avec l'inverse de $f$.\end{demo}

\begin{remarque3}
Le diagramme $(\Xsi^h \longrightarrow \Xsib^h \longleftarrow Z_\sigma)$ ne définit pas un objet de $\mathsf{OuvComp}_\mathcal{U}$ car $\Xsib^h$ n'est pas complet.
\end{remarque3}

\subsection{Cartes locales algébriques}

Le schéma affine $\Xsib^h$ est limite projective de schémas étales sur $\mathcal{U}\times_{\mathcal{B}_\sigma} \Xsib$ et le schéma semi-abélien $G_\sigma^h$ est de type fini sur $\Xsib^h$. Il existe donc un schéma $\Xsib^{et}$, un morphisme étale
$$\Xsib^{et} \longrightarrow \mathcal{U}\times_{\mathcal{B}_\sigma} \Xsib$$
qui est un isomorphisme sur $Z_\sigma$, et un schéma semi-abélien $G_\sigma^{et}$ sur $\Xsib^{et}$ qui descend $G_\sigma^h$.
On pose $$\Xsi^{et} = \Xsi \times_{\Xsib} \Xsib^{et}\: \: \mathrm{et} \: \: \partial \Xsib^{et} = \Xsib^{et} - \Xsi^{et}\: .$$
Le schéma semi-abélien $G_\sigma^{et}$ vérifie les mêmes propriétés que celles de  $G_\sigma^h$  données dans la proposition~\ref{propGxih}. Par exemple, il est abélien principalement polarisé sur $\Xsi^{et}$ et son groupe de $p$-torsion est canoniquement isomorphe à $$M_\sigma[p]\times_{\Xsi} \Xsi^{et}\: .$$
Si $(U\rightarrow S \leftarrow S_0)$ est un objet de $\mathsf{OuvComp}_{\mathcal{U}}$, toute flèche
$$(U\longrightarrow S \longleftarrow S_0) \longrightarrow (\Xsi^{et} \longrightarrow \Xsib^{et} \longleftarrow Z_\sigma )$$
se factorise par $\Xsib^h$ car $(S,S_0)$ est hensélien. On déduit de la proposition~\ref{propModProchFrontHens} une description des points proches de la frontière de $$(\Xsi^{et} \longrightarrow \Xsib^{et} \longleftarrow Z_\sigma )\: .$$
Le schéma $\Xsib^{et}$ n'est pas unique mais les compactifications ne dépendront pas du choix.
Le schéma abélien principalement polarisé $G_{\sigma}^{et}\rightarrow \Xsi^{et}$ induit une flèche
$\Xsi^{et} \rightarrow \ag \: .$
Appelons \textit{carte locale algébrique} associée à  $\sigma\in \Sigma$ tout diagramme  $$\left(\Xsib^{et} \longleftarrow \Xsi^{et} \longrightarrow  \ag \right)$$
obtenu par la construction précédente.

\subsubsection{Cartes locales algébriques avec niveau parahorique}

Notons $\XsiKb^{et}$ le produit fibré 
\begin{eqnarray*}
\xymatrix{\XsiKb^{et} \ar[r] \ar[d] & \Xsib^{et} \ar[d] \\
\XsiKb \ar[r] & \Xsib\: .}
\end{eqnarray*}
Il contient comme sous-schéma ouvert $\XsiK^{et}$, qui est défini par  le produit fibré
\begin{eqnarray*}
\xymatrix{\XsiK^{et} \ar[r] \ar[d] & \Xsi^{et} \ar[d] \\
\XsiK \ar[r] & \Xsi\: .}
\end{eqnarray*}
La proposition suivante résulte de l'approximation préservant le groupe de $p$-torsion.

\begin{proposition3} \label{propXsiKXsiEspMod} Le morphisme $\XsiK^{et} \rightarrow \Xsi^{et}$ est l'espace de modules relatif des structures de niveau parahoriques de type $\mathcal{D}$ en $p$ sur $G_{\sigma}^{et}[p]$.
\end{proposition3}

\begin{demo} Le morphisme $\XsiK\rightarrow \Xsi$ est l'espace de modules relatif des structures de niveau parahoriques sur $M_{\sigma}[p]$, et le schéma abélien $G_{\sigma}^{et} \rightarrow \Xsi^{et}$ est muni d'un isomorphisme $G_{\sigma}^{et}[p] \iso M_{\sigma}[p] \times_{\Xsi} \Xsi^{et}$.
\end{demo}

Nous avons obtenu un diagramme à carrés cartésiens et flèches verticales propres
\begin{eqnarray*}
\xymatrix{ \XsiKb^{et} \ar[d] & \XsiK^{et} \ar[l] \ar[r] \ar[d] & \agK \ar[d]\\
\Xsib^{et} & \Xsi^{et} \ar[l] \ar[r] & \ag\: .}
\end{eqnarray*}
On appelle ce diagramme  \textit{carte locale algébrique de niveau parahorique de type $\mathcal{D}$ en $p$} associée à $\sigma$.  Il correspond une unique carte locale algébrique avec niveau à toute carte locale algébrique sans niveau.

\begin{remarque3} Nous n'avons pas besoin de définir des cartes locales formelles avec niveau et de les approximer. La mauvaise réduction empêcherait de toute façon le contrôle d'une éventuelle approximation.
\end{remarque3}

\section{Compactifications toroïdales}

\subsection{Recollement des cartes locales algébriques}

On choisit pour tout $\sigma \in \Sigma$  une carte locale algébrique $$(\xymatrix{\Xsib^{et} & \Xsi^{et} \ar[l]\ar[r] & \ag})$$ et l'on considère
$$\overline{U} = \coprod_{\sigma\in\Sigma} \Xsib^{et},$$
qui est localement de type fini sur $\Spec(\Z)$. On pose 
$${U} = \coprod_{\sigma\in\Sigma} \Xsi^{et}\: , \: Z=\coprod_{\sigma\in\Sigma} Z_{\sigma} \: \:
\mathrm{et} \:\: \partial \overline{U} = \overline{U} - U\: , $$
qui est un diviseur à croisements normaux de $\overline{U}$. On note $\underline{X}=\coprod \underline{X}_\sigma$, qui est un faisceau constructible sur $U$, et $\underline{B}=\coprod \underline{B}_\sigma$, qui est une forme bilinéaire sur $\underline{X}$ à valeurs dans le faisceau $\underline{\mathrm{Div}}_{\overline{U}}$.

On a construit un schéma semi-abélien $G=\coprod G_\sigma$ sur $\overline{U}$ qui est abélien principalement polarisé sur $U$. Il existe des isomorphismes $\underline{X}(G)\iso \underline{X}$ et $\underline{B}(G)\iso \underline{B}$, et l'application de Kodaira-Spencer induit un isomorphisme
$$K_{G/\overline{U}}\: : \: \mathrm{Sym}^2 (\Omega_G) \isolong \Omega^1_{\overline{U}/\Z}\left(\mathrm{log} \: \partial \overline{U}\right) \: .$$
Le schéma abélien $G$ sur $U$ induit un diagramme
$$\left(\xymatrix{\overline{U} & U \ar[l] \ar[r] & \ag }\right) \: .$$
La flèche $U\rightarrow \ag$ est surjective car $\mathcal{M}_{\{0\}}=\ag$.

\begin{proposition2} \label{propUagEt} Le morphisme $U\rightarrow \ag$ est étale.
\end{proposition2}

\begin{demo} On sait que le schéma abélien universel $\Ag\rightarrow \ag$ induit un isomorphisme de Kodaira-Spencer
$$K_{\Ag/\ag} \: : \: \mathrm{Sym}^2(\Omega_{\Ag})\isolong \Omega^1_{\ag/\Z}$$
et le schéma abélien $G\rightarrow U$ un isomorphisme
$$K_{G/U} \: : \: \mathrm{Sym}^2 (\Omega_G)\: \isolong \:\Omega^1_{U/\Z}\:. $$
Notons $f:U\rightarrow \ag$ le morphisme induit par $G$. On a $G=f^* \Ag$ et l'on a $f^* \: \mathrm{Sym}^2(\Omega_{\Ag}) \:\simeq \:\mathrm{Sym}^2 (\Omega_G)$
et $f^* K_{\Ag/\ag} = K_{G/U}$ donc
$$f^*\: \Omega^1_{\ag / \Spec(\Z)}\: \isolong \:\Omega^1_{U/\Spec(\Z)}\: .$$ 
Comme les champs algébriques $U$ et $\ag$ sont lisses sur $\Spec(\Z)$,  le critère jacobien montre que $f$ est étale.
\end{demo}

Considérons le schéma $$R=U\times_{\ag} U=\underline{\mathrm{Isom}}_{U\times U}(\mathrm{pr}_1^*\: G,\mathrm{pr}_2^*\:G)\: .$$ Il est muni de deux flèches étales vers $U$, et le champ algébrique $\ag$ est isomorphe au champ quotient $[U/R]$. Notons $\overline{R}$ le normalisé de $$\overline{U}\times_{\Spec(\Z)} \overline{U}$$ dans $R$. C'est un schéma muni de deux flèches $\overline{R}\rightrightarrows \overline{U}$ prolongeant $R\rightrightarrows U$. Il existe un isomorphisme $\mathrm{pr}_1^*\: G \iso \mathrm{pr}_2^*\: G$ au-dessus de $R$. D'après le théorème~\ref{thRay}, il s'étend en un isomorphisme au-dessus de $\overline{R}$ par normalité de ce schéma. La démonstration de la proposition suivante est tirée de~\cite{Deg@FaltingsChai} et servira de modèle à celle de la proposition~\ref{propR0bet}.

\begin{proposition2}[\cite{Deg@FaltingsChai} IV.5.4] \label{propRbEt}  Les deux flèches $\overline{R}\rightrightarrows \overline{U}$ sont étales.
\end{proposition2}

\begin{demo} Soient $z$ un point de $\overline{R}$ et $x$ et $y$ ses images dans $\overline{U}$. Il suffit de montrer que les complétés des localisés de $\overline{R}$ en $z$ et de $\overline{U}$ en $x$ et $y$ sont isomorphes \textit{via} les deux flèches. Comme les flèches $R \rightrightarrows U$ sont étales, on peut supposer que $x$ et $y$ sont des éléments du bord $\partial\overline{U}$. Il existe  des cônes $\sigma_x, \: \sigma_y \in \Sigma$ tels que $x$ soit sur la $\sigma_x$-strate et $y$ sur la $\sigma_y$-strate. Soient $\mathcal{O}_x$, $\mathcal{O}_y$ et $\mathcal{O}_z$ les anneaux des complétés des localisés en $x$, $y$ et $z$. Les morphismes $\overline{R} \rightrightarrows \overline{U}$ sont dominants car les morphismes ${R} \rightrightarrows {U}$ sont surjectifs. On dispose d'injections d'anneaux
$$\xymatrix{ & \mathcal{O}_z & \\
\mathcal{O}_x \ar@{_{(}->}[ru] & & \mathcal{O}_y \ar@{^{(}->}[lu]}$$
et il s'agit de montrer que ces injections sont des égalités.

On note $G_x$ et $G_y$ les schémas semi-abéliens sur $\Spec(\mathcal{O}_x)$ et $\Spec(\mathcal{O}_y)$ obtenus par image inverse de $G$. D'après la proposition~\ref{propmodfrontXsib}, le $1$-motif de Mumford associé à $G_x$ induit une flèche $\Spec(\mathcal{O}_x)\rightarrow {{\overline{\mathcal{M}}}_{\sigma_x}}$ et celui associé à $G_y$ une flèche $\Spec(\mathcal{O}_y)\rightarrow {\overline{\mathcal{M}}}_{\sigma_y}$.
Le cône $\sigma_x$ est déterminé par $\underline{B}(G_x)$ car c'est l'unique cône minimal parmi ceux dont l'adhérence contient les $p_D\circ \underline{B}(G_x)$, où $D$ est un diviseur irréductible de $\Spec(\mathcal{O}_x)$. Il en est de même pour $\sigma_y$ et comme $\underline{B}(G_x)|_{\mathcal{O}_z} \simeq \underline{B}(G_y)|_{\mathcal{O}_z},$ les cônes $\sigma_x$ et $\sigma_y$ sont conjugués sous l'action de $\mathrm{GL}(X)$ et il existe un isomorphisme $\overline{\mathcal{M}}_{\sigma_x} \iso \overline{\mathcal{M}}_{\sigma_y}$. On obtient un diagramme commutatif
$$\xymatrix{ & \Spec(\mathcal{O}_x) \ar[r] & {{\overline{\mathcal{M}}}_{\sigma_x}} \ar[dd]_{\sim} \\
\Spec(\mathcal{O}_z) \ar[ru] \ar[rd] & & \\
& \Spec(\mathcal{O}_y) \ar[r] & {\overline{\mathcal{M}}}_{\sigma_y}  }$$
Les anneaux $\mathcal{O}_x$ et $\mathcal{O}_y$ sont les complétés des localisés de ${{\overline{\mathcal{M}}}_{\sigma_x}}$ et ${{\overline{\mathcal{M}}}_{\sigma_y}}$ en les images du point fermé de $\Spec(\mathcal{O}_z)$, donc il existe un isomorphisme
$$\phi:\Spec(\mathcal{O}_x)\isolong \Spec(\mathcal{O}_y)\: .$$
L'action de $\mathrm{GL}(X)$ induit un isomorphisme entre les $1$-motifs de Mumford universels sur ${\overline{\mathcal{M}}}_{\sigma_x}$ et ${\overline{\mathcal{M}}}_{\sigma_y}$, et il existe un isomorphisme $G_x \iso \phi^*G_y$ sur $\Spec(\mathcal{O}_x)$ qui descend l'isomorphisme $G_x\iso G_y$ sur $\Spec(\mathcal{O}_z)$.

Pour démontrer la proposition, on construit une rétraction de l'injection $\mathcal{O}_x \hookrightarrow \mathcal{O}_z$. L'isomorphisme $G_x \iso \phi^* G_y$ sur $\Spec(\mathcal{O}_x)$ induit un diagramme commutatif
$$\xymatrix{R\ar[r]\ar[d] & \overline{R} \ar[d] \\
U\times U \ar[r] & \overline{U}\times \overline{U} \\
\Spec(\mathrm{Frac}(\mathcal{O}_x)) \ar@/^2pc/[uu] \ar[r] \ar[u] & \Spec(\mathcal{O}_x) \ar[u]_{(\mathrm{Id},\phi)}}$$
L'anneau $R$ est normal et la flèche $\Spec(\mathrm{Frac}(\mathcal{O}_x))\rightarrow R$ s'étend en un morphisme $\Spec(\mathcal{O}_x)\rightarrow \overline{R}.$
Cette flèche se factorise par le complété de $\overline{R}$ en $z$ et induit bien une  rétraction de $\mathcal{O}_x \hookrightarrow \mathcal{O}_z$.
\end{demo}

\begin{remarque2} \label{remRbet} Nous avons utilisé pour la première fois le fait que les cônes de $\Sigma$ ne se recouvrent pas.
\end{remarque2}

\begin{lemme2}[\cite{Deg@FaltingsChai} corollaire IV.5.5]  La double-flèche $\overline{R}\rightrightarrows \overline{U}$ est un groupoïde étale.
\end{lemme2}

\begin{demo} Il suffit de trouver une flèche $\overline{R}\times_{\overline{U}} \overline{R} \rightarrow \overline{R}$ qui étend la flèche $R\times_U R\rightarrow R$ et qui vérifie les axiomes d'un groupoïde. Comme le schéma $\overline{R}\times_{\overline{U}} \overline{R}$ est étale sur $\overline{U}$, il est égal à la normalisation de $\overline{R}$ dans $R\times_U R$. La flèche $\overline{R}\times_{\overline{U}} \overline{R} \rightarrow \overline{R}$ provient de la propriété universelle de la normalisation et les axiomes sont vérifiés par densité de $U$ dans $\overline{U}$.
\end{demo}

La flèche $R\rightarrow U\times U$ est finie car $\ag$ est séparé, et le morphisme $\overline{R}\rightarrow \overline{U}\times \overline{U}$ est également fini car $\overline{U}$ est excellent.

\subsubsection{Description des compactifications sans niveau}

Considérons le champ algébrique quotient $${\overline{\mathcal{A}}_g^{\Sigma}} = [\overline{U}/\overline{R}]\: .$$
Il est séparé, lisse sur $\Spec(\Z)$, contient $\ag$ comme ouvert dense et $$\partial {\overline{\mathcal{A}}_g^{\Sigma}} \: = \: {\overline{\mathcal{A}}_g^{\Sigma}} - \ag$$ comme diviseur à croisements normaux. Le schéma semi-abélien $G\rightarrow \overline{U}$ est muni d'une donnée de descente $\mathrm{pr}_1^* G \iso \mathrm{pr}_2^* G$ sur $\overline{R}$ donc induit un schéma semi-abélien  $$G\longrightarrow {\overline{\mathcal{A}}_g^{\Sigma}}$$ qui étend la variété abélienne universelle. L'application de Kodaira-Spencer induit un isomorphisme
$$\mathrm{Sym}^2 (\Omega_G)\: \isolong \:\Omega^1_{{\overline{\mathcal{A}}_g^{\Sigma}}/\Z}(\mathrm{log}\: \partial {\overline{\mathcal{A}}_g^{\Sigma}})\:. $$
La stratification  induite  sur  ${\overline{\mathcal{A}}_g^{\Sigma}}$ par les branches du diviseur à croisements normaux est paramétrée par l'ensemble fini $\Sigma/\mathrm{GL}(X)$. La complétion formelle 
le long de la strate associée $\sigma\in \Sigma$ est isomorphe au quotient de la complétion formelle de $\Xsib$ le long de $Z_\sigma$ par le groupe fini
$$\Gamma_\sigma \: = \: \mathrm{Stab}_{\mathrm{GL}(X_\sigma)} (\sigma) \: .$$
Le schéma semi-abélien $G$ correspond par cet isomorphisme au quotient de Mumford $G_\sigma$ du $1$-motif $M_\sigma\rightarrow \Xsi$, quitte à composer par l'automorphisme de recomplétion de la proposition~\ref{propRecomplAutoEt}. On renvoie le lecteur au théorème~IV.5.7 de~\cite{Deg@FaltingsChai} pour plus détails.
Soit $(V\hookrightarrow S)$ un objet de $\mathsf{Ouv}$. L'image inverse de $G$ définit un foncteur de $(\ag\hookrightarrow {\overline{\mathcal{A}}_g^{\Sigma}})(V\hookrightarrow S)$ vers la catégorie des schémas semi-abéliens sur $S$.

\begin{proposition3} \label{propModFrontAgb} Ce foncteur réalise une équivalence de catégories entre 
$$\left(\ag\hookrightarrow {\overline{\mathcal{A}}_g^{\Sigma}}\right)(U\hookrightarrow S)$$
et la catégorie des schémas semi-abéliens $G'\rightarrow S$  qui sont abéliens principalement polarisés de genre $g$ sur $U$, et tels que localement pour la topologie étale de $S$, le faisceau $\underline{X}(G')$ est quotient de $X$ et il existe $\sigma \in \Sigma$ dont l'adhérence contient la forme bilinéaire $p_D\circ\underline{B}(G')$ pour tout diviseur irréductible réduit $D$ de $S$.
\end{proposition3}

\begin{demo} Le schéma abélien $G'$ principalement polarisé sur $V$ induit une flèche $V\rightarrow \ag$ que l'on veut prolonger. Un éventuel prolongement est unique car $\overline{\ag}$ est séparé, et son existence se teste localement. On peut donc supposer que $S$ est le spectre d'un anneau local complet de point fermé $S_0$. On choisit $\sigma\in \Sigma$ minimal parmi les cônes dont l'adhérence contient les $\underline{B}(G')\otimes p_D$
et l'on relève la flèche $V\rightarrow \ag$ en une flèche $V\rightarrow \Xsi^{et}$ qui s'étend en $$S\longrightarrow \Xsib^{et}$$ grâce à la propriété modulaire proche de la frontière de $\Xsib^{et}$.
\end{demo}

On en déduit que le champ ${\overline{\mathcal{A}}_g^{\Sigma}}$ ne dépend ni du choix de la présentation $\mathcal{U}$ de $\mathcal{B}_\sigma$ ni du choix des approximations $G^{et}_\sigma$. 
Il reste à montrer la proposition suivante pour vérifier que l'on a construit une compactification de $\ag$.

\begin{proposition3} \label{propPROPRE} Le champ algébrique ${\overline{\mathcal{A}}_g^{\Sigma}}$ est propre sur $\Spec(\Z)$.
\end{proposition3}

\begin{demo} La flèche est séparée de type fini et l'on peut appliquer le critère valuatif de propreté de~\ega{ii}{7.3.10}. Soient $S$ un trait de point générique $\eta$ et $\eta \rightarrow \ag$ un morphisme. On veut le prolonger en une flèche
$$S \longrightarrow {\overline{\mathcal{A}}_g^{\Sigma}}$$
localement pour la topologie \textit{\textit{fppf}}. Soit $G'$ la variété abélienne sur $\eta$ induite par $\eta\rightarrow \ag$. Le théorème de réduction semi-stable pour les variétés abéliennes assure l'existence d'un schéma semi-abélien $G'\rightarrow S$ prolongeant $G'\rightarrow \eta$, quitte à effectuer un changement ramifié de trait. Un changement étale de trait permet de supposer que $\underline{X}(G)$ est quotient de $X$. Un trait admet un unique diviseur effectif réduit irréductible, la flèche $\underline{B}(G) :  X\otimes X \rightarrow \Z$ est semi-définie positive à radical rationnel d'après la proposition \ref{propBGsab} et $\Sigma$ recouvre le cône  des formes semi-définies positives à radical rationnel. La propriété modulaire à la frontière de 
$$\left(\ag \hookrightarrow {\overline{\mathcal{A}}_g^{\Sigma}}\right)$$
montre alors l'existence d'un prolongement.
\end{demo}

\begin{remarque3} Le contrôle de l'approximation a été indispensable pour vérifier la propreté de ${\overline{\mathcal{A}}_g^{\Sigma}}$.
\end{remarque3}

\subsubsection{Cas d'une structure de niveau}

Dans ce paragraphe, on définit la compactification avec niveau parahorique associée à $\Sigma$. On associe à chaque carte locale algébrique sans niveau $$\xymatrix{\Xsib^{et} & \Xsi^{et} \ar[l]\ar[r] & \ag}$$ la carte locale algébrique de niveau parahorique de type $\mathcal{D}$ en $p$  $$\xymatrix{ \XsiKb^{et} \ar[d] & \XsiK^{et} \ar[l] \ar[r] \ar[d] & \agK \ar[d]\\
\Xsib^{et} & \Xsi^{et} \ar[l] \ar[r] & \ag}$$qui lui correspond. Considérons les schémas 
$$\overline{U}_0 \:=\: \coprod_{\sigma\in\Sigma} \XsiKb^{et} \quad \mathrm{et} \quad {U}_0 \:=\: \coprod_{\sigma\in\Sigma} \XsiK^{et}\: .$$
Ils induisent un diagramme cartésien à flèches verticales propres
\begin{eqnarray*}
\xymatrix{ \overline{U}_0 \ar[d] & U_0 \ar[l] \ar[r] \ar[d] & \agK \ar[d]\\
\overline{U} & U\ar[l] \ar[r] & \ag\: ,}
\end{eqnarray*}
par construction des cartes locales algébriques avec niveau. Définissons le groupoïde $R_0 \rightrightarrows U_0$ comme le produit fibré $$R_0=U_0 \times_{\agK} U_0\: .$$
Les deux flèches $R_0 \rightrightarrows U_0$ sont étales car $U\rightarrow \ag$ est étale, et le champ quotient $[U_0/R_0]$ est bien sûr isomorphe à $\agK$.

\begin{remarque3} Pour démontrer que $U_0\rightarrow \agK$ est étale, on n'a pas procédé directement. Du fait de la mauvaise réduction, il semble difficile d'adapter la technique de~\cite{Deg@FaltingsChai}, basée sur l'isomorphisme de Kodaira-Spencer.
\end{remarque3}

\subsubsection{Normalité} Dans ce paragraphe, nous démontrons la normalité de $\overline{\mathcal{M}}_{\sigma,\: 0}$, qui sera utile pour la proposition~\ref{propR0bet}. Nous utilisons les résultats de~\cite{FlatSympl@Gortz}, selon lesquels le morphisme $\agK \longrightarrow \Spec(\Z_p)$ est plat et sa fibre spéciale est réduite.

\begin{proposition3} \label{propAgKNorm} Le champ algébrique $\agK$ est normal.
\end{proposition3}

\begin{demo} D'après le critère de normalité de Serre, il suffit de montrer que $\agK$ est régulier en codimension $1$ et que ses points de codimension~$\geq 2$ sont de profondeur~$\geq 2$.
D'après~le théorème principal de~\cite{VarAb@GenestierNgo}, le lieu régulier de $\agK$ coïncide avec le lieu ordinaire, et son complémentaire dans $\agK$ est de codimension~$\geq 2$. Cela montre que $\agK$ est régulier en codimension~$1$.

Vérifions la deuxième condition du critère de normalité de Serre. Soient $A$ l'anneau local de $\agK$ en un point $x$ de profondeur $\geq 2$ et $\mathfrak{m}$ son idéal maximal. Comme $\agK\rightarrow \Spec(\Z[1/p])$ est lisse, on peut supposer que $p\in \mathfrak{m}$. L'anneau $A$ est plat sur $\Z_{(p)}$ donc $p$ n'est pas diviseur de~$0$ dans~$A$. On note $\overline{A}=A/pA$ et $\overline{\mathfrak{m}}=\mathfrak{m}/\mathfrak{m}\cap pA$. D'après~\cite{FlatSympl@Gortz}, $\overline{A}$ est réduit donc sans composante immergée. Les idéaux associés de $\overline{A}$ sont donc minimaux. On en déduit que $\overline{\mathfrak{m}}$  n'est pas associé. D'après la proposition~8 du paragraphe~IV.1.4 de~\cite{AC34@Bourbaki}, il existe $x\in \overline{\mathfrak{m}}$ qui n'est pas diviseur de $0$ dans $\overline{A}$. La suite de paramètres $(p,x)$ est bien de longueur $2$.
\end{demo}

Soit $\sigma \in \Sigma$. Le corollaire suivant résulte de la proposition~\ref{propAgKNorm}, des corollaires~\ref{coroAb} et~\ref{corollaireLibre}, et de la normalité des plongements toriques.

\begin{corollaire3} \label{coroXsiKbNormal} Le champ algébrique $\XsiKb$ est normal.
\end{corollaire3}

\subsubsection{Relation d'équivalence avec niveau} \label{partieRelEqNiveau} Dans ce paragraphe, on construit une relation d'équivalence étale sur $\overline{U}_0$.
On note $\overline{R}_0$ le normalisé de $\overline{U}_0 \times \overline{U}_0$ dans $R_0$. Il est muni de deux flèches 
$$\overline{R}_0\rightrightarrows \overline{U}_0$$ qui étendent $$R_0\rightrightarrows U_0\: .$$

La restriction du schéma semi-abélien $G$ à $U_0$ est abélien, muni d'un drapeau $H_{\bullet}\subset G[p]$. On dispose sur $R_0$ d'un isomorphisme
$$\mathrm{pr}_1^* \: G \: \isolong \:  \mathrm{pr}_2^* \: G$$
qui envoie $\mathrm{pr}_1^* \: H_{\bullet}$ sur $\mathrm{pr}_2^* \: H_{\bullet}$. Par normalité, cet isomorphisme se prolonge en un isomorphisme sur $\overline{R}_0$.

\begin{proposition3} \label{propR0bet} Les deux flèches $\overline{R}_0\rightrightarrows \overline{U}_0$ sont étales.
\end{proposition3}

\begin{demo} Soit $z$ un point de $\overline{R}_0$ et $x$, $y$ ses deux images dans $\overline{U}_0$. Montrons que les complétés $\mathcal{O}_z$, $\mathcal{O}_x$ et $\mathcal{O}_y$ des localisés de $\overline{R}_0$ en $z$ et de $\overline{U}_0$ en $x$ et $y$ sont isomorphes. On peut supposer que $x$ et $y$ sont des éléments du bord. Il existe des cônes $\sigma_x$ et $\sigma_y$ de $\Sigma$ tels que $x$ appartienne à la $\sigma_x$-strate et $y$ sur la $\sigma_y$-strate. On dispose d'injections d'anneaux
$$\xymatrix{ & \mathcal{O}_z & \\
\mathcal{O}_x \ar@{_{(}->}[ru] & & \mathcal{O}_y \ar@{^{(}->}[lu]}$$
et de schémas semi-abéliens $G_x$ et $G_y$ obtenus par image inverse de $G$ sur $\Spec(\mathcal{O}_x)$ et  $\Spec(\mathcal{O}_x)$. Ils sont génériquement abéliens, munis d'une structure de niveau parahorique $H_{\bullet,\: x}$ et $H_{\bullet,\: y}$, et il existe sur $\Spec(\mathcal{O}_z)$ un isomorphisme
$$\mathrm{pr}_1^*\: G_x  \: \isolong \: \mathrm{pr}_2^* \: G_y$$
qui respecte les structures de niveau.
Les anneaux $\mathcal{O}_x$ et $\mathcal{O}_y$ sont normaux d'après le corollaire~\ref{coroXsiKbNormal}. D'après la proposition~\ref{propMumIw}, on obtient deux $1$-motifs de Mumford avec structure de niveau associés à $G_x$ et $G_y$, et la proposition~\ref{propmodfrontXsiKb} montre l'existence de flèches $$\Spec(\mathcal{O}_x)\longrightarrow {\overline{\mathcal{M}}}_{\sigma_x,\: 0} \: \: \mathrm{et} \: \:\Spec(\mathcal{O}_y)\longrightarrow {\overline{\mathcal{M}}}_{\sigma_x,\: 0}$$
par lesquelles descendent les deux $1$-motifs avec structure de niveau.
On raisonne exactement comme dans la preuve de la proposition~\ref{propRbEt} pour construire une rétraction de l'injection $\mathcal{O}_x\hookrightarrow \mathcal{O}_z$.
\end{demo}

\medskip

On voit comme précédemment que $\overline{R}_0\rightrightarrows \overline{U}_0$ est un groupoïde étale. Appelons
$${\overline{\mathcal{A}}_{g,\: 0}^{\Sigma}}\: =\: [\: \overline{U}_0/\overline{R}_0 \:]$$
le champ algébrique quotient. Expliquons rapidement la propriété modulaire à la frontière vérifiée par $$\left( \agK \: \hookrightarrow\: {\overline{\mathcal{A}}_{g,\: 0}^{\Sigma}} \right) \: .$$
Soient $(U\hookrightarrow S)$ un objet de $\mathsf{Ouv}$ et $(G',H'_\bullet)$ un objet de $\agK(U)$. On suppose que $G'$ s'étend en un schéma semi-abélien sur $S$. Le quotient $G'/H'_i$ est alors représentable par un schéma semi-abélien sur $S$ pour $0\leq i  \leq s$. Supposons que $\underline{X}(G')$ soit un quotient de $X$. D'après la proposition~\ref{propXi0w} et le diagramme~(\ref{equationFlip}), la famille des formes bilinéaires $\underline{B}(G'/H'_i)$ définit un morphisme $\underline{B}(G',H'_\bullet)$ de $\mathrm{Sym}^2(X\otimes \Q)$ dans $\underline{\mathrm{Div}}_S\otimes \Q$. La propriété modulaire à la frontière suivante
se montre comme la proposition~\ref{propModFrontAgb}, en utilisant la proposition~\ref{propmodfrontXsiKb}.

\begin{proposition3} \label{propModFrontAgKbVersion1} Le groupoïde $\left(\agK\hookrightarrow {\overline{\mathcal{A}}_{g,\: 0}^{\Sigma}}\right)(U\hookrightarrow S)$
paramètre les objets $(G',H'_\bullet)$ de $\agK(U)$ tels que $G'$ s'étend en un schéma semi-abélien sur $S$, et tels que localement pour la topologie étale de $S$, le faisceau $\underline{X}(G)$ est quotient de $X$ et il existe $\sigma \in \Sigma$ dont l'adhérence contient la forme bilinéaire $p_D \circ \underline{B}(G',H'_\bullet)$ pour tout diviseur irréductible réduit $D$ de $S$.
\end{proposition3}

Le théorème suivant résume les premières propriétés de notre construction.

\begin{theoreme3} \label{thPrinc1} Le champ algébrique ${\overline{\mathcal{A}}_{g,\: 0}^{\Sigma}}$ ne dépend que de $\Sigma$. Il est séparé,  normal et sa projection vers $\Spec(\Z)$ est de type fini, plate, à fibres géométriquement réduites.
Ce champ contient $\agK$ comme ouvert dense et il existe un diagramme cartésien 
$$\xymatrix{ {\agK} \ar[r]\ar[d] & {\overline{\mathcal{A}}_{g,\: 0}^{\Sigma}} \ar[d]\\
\ag \ar[r] & {\overline{\mathcal{A}}_g^{\Sigma}} }$$
à flèches verticales propres. En particulier, ${\overline{\mathcal{A}}_{g,\: 0}^{\Sigma}}$ est propre sur $\Spec(\Z)$ et est une compactification de $\agK$.
Il existe un schéma semi-abélien $G$ et un drapeau de groupes quasi-finis plats $$H_\bullet \subset G[p]$$  sur ${\overline{\mathcal{A}}_{g,\: 0}^{\Sigma}}$ qui étendent le schéma abélien et le drapeau universel sur $\agK$.
\end{theoreme3}

\begin{demo} Les assertions sont des conséquences immédiates de la construction. La proposition~\ref{propModFrontAgKbVersion1} montre que les compactifications ne dépendent que de $\Sigma$. La platitude et la réduction des fibres des compactifications proviennent des résultats de~\cite{FlatSympl@Gortz}. Le drapeau s'obtient par adhérence du drapeau universel sur $\agK$ dans le groupe quasi-fini plat $G[p]$.
\end{demo}

Le champ ${\overline{\mathcal{A}}_{g,\: 0}^{\Sigma}}$ est muni d'une stratification paramétrée par un quotient fini de $\Sigma$ et la complétion formelle le long de la strate associée à $\sigma\in \Sigma$ est isomorphe au quotient de la complétion formelle de $\XsiKb$ le long de $Z_{\sigma,\: 0}$ par un groupe fini. Le schéma semi-abélien $G$ et le drapeau $H_\bullet$ correspondent par cet isomorphisme au quotient de Mumford du $1$-motif $M_\sigma$ muni de sa structure de niveau  canonique sur $\XsiK$. Dans le paragraphe suivant, on précise la combinatoire du bord de la compactification.

\subsection{Compactifications améliorées} \label{partieCompactAmeliore} Deux problèmes subsistent. Le premier est de décrire précisément le bord de $${\overline{\mathcal{A}}_{g,\: 0}^{\Sigma}}$$
et notamment les relations d'incidence entre les différentes strates. Le second est de se débarrasser des singularités toriques intervenant dans la compactification pour obtenir un champ n'ayant que les singularités des $\agrK$ pour $r\leq g$. Ce champ sera lisse sur $\Spec(\Z[1/p])$ et le complémentaire de $\ag$ y sera un diviseur à croisements normaux.

Ces problèmes se résolvent simultanément en introduisant un \textit{complexe conique} $\mathcal{C}$ contenant $C(X)$. La décomposition $\Sigma$ de $C(X)$ induit une décomposition $\mathfrak{S}_\Sigma$ de~$\mathcal{C}$, ce qui conduit à une description explicite du bord de $${\overline{\mathcal{A}}_{g,\: 0}^{\Sigma}}\: .$$
On peut également remplacer le choix de $\Sigma$ par un choix plus fin d'une décomposition~$\mathfrak{S}$ de~$\mathcal{C}$, ce qui permet de construire des compactifications améliorées, sans singularité d'origine torique.

\subsubsection{Sous-espaces isotropes} On introduit des notations analogues à celles du paragraphe~\ref{partiePosRel}. Soit $V=\oplus_{j=1}^{2g}\Z \cdot x_j$ un~$\Z$-module libre de rang~$2g$  muni d'une base~$(x_j)$. On le munit de la forme alternée donnée par la matrice
$$J=\left( \begin{array}{cc} 0 & J' \\
 -J' & 0           \end{array} \right) \: ,$$
où~$J'$ est la matrice  anti-diagonale dont  les coefficients non nuls sont égaux à~$1$. 
On rappelle que $\mathcal{D}=\lbrace d_1 < \cdots < d_s \rbrace$.
On définit une chaîne parahorique autoduale $\mathbb{V}_\mathcal{D}^\bullet$ de $V$
en posant pour $1\leq i\leq s$,
$$\mathbb{V}^i_\mathcal{D}  \:= \: \bigoplus_{j=1}^{d_i}\: \Z \cdot x_j \: \oplus \bigoplus_{j=d_i+1}^{2g}p\: \Z \cdot x_j$$
$$\mathrm{et} \:\:\: \: \: \mathbb{V}^{s+i}_\mathcal{D}\: = \: \bigoplus_{j=1}^{2g-d_{s+1-i}} \:  \Z\cdot x_j \: \oplus \bigoplus_{j=2g - d_{s+1-i}+1}^{2g} p\:\Z\: .$$
Notons $\Gamma_0$  le stabilisateur de $\mathbb{V}_\mathcal{D}^\bullet$ dans $\gsp(\Z)$ : c'est le sous-groupe de congruence parahorique associé à~$\mathcal{D}$.
Si $V'$ est un sous-espace totalement isotrope facteur direct de $V$, notons $r$ son rang et $V'^\perp$ son orthogonal. L'entier $r$ détermine l'ensemble de positions relatives $\mathcal{W}_r$ (\textit{cf.} le paragraphe~\ref{partiePosRel}, où cet ensemble était noté~$\mathcal{W}$).
L'image de $\mathbb{V}_\mathcal{D}^\bullet$ dans $V/V'^\perp$ définit un diagramme commutatif
$$\xymatrix{ \mathbb{Y}_0 \ar[r] \ar[d] & \mathbb{Y}_1 \ar[r] \ar[d] & \cdots  \ar[r] & \mathbb{Y}_s \ar[d] \\
\mathbb{X}_0 \ar[d]  & \mathbb{X}_1 \ar[l]\ar[d] & \cdots \ar[l] & \mathbb{X}_s\ar[d] \ar[l] \\
\mathbb{Y}_0 \ar[r] & \mathbb{Y}_1 \ar[r] & \cdots \ar[r]  & \mathbb{Y}_s }$$
En particulier, $V'$ définit un élément $w\in \mathcal{W}_r$. On est dans le cadre d'application de la construction du paragraphe~\ref{partieVariante}, et l'on peut considérer le champ de modules 
$$\mathcal{M}_{V',\: 0} = \mathcal{M}_{\mathbb{X}_\bullet,\mathbb{Y}_\bullet,\: 0}^w\: .$$
D'après les résultats du paragraphe~\ref{partieEspModPar}, il existe un champ algébrique $\mathcal{B}_{V',\: 0}$ et un tore $E_{V',\: 0}$ tels que $\mathcal{M}_{V',\: 0}$ soit un $E_{V',\: 0}$-torseur sur $\mathcal{B}_{V',\: 0}$. Posons $S_{V',\: 0}=\Hom(E_{V',\: 0},\Gm)$. C'est une structure entière sur $\mathrm{Sym}^2(V/V'^\perp)\otimes \R$.

\subsubsection{Complexe conique} Notons $\mathfrak{C}$ l'ensemble des sous-espaces totalement isotropes facteurs directs de $V$. Pour tout $V'\in \mathfrak{C}$, notons  $C(V/V'^\perp)$ le cône des formes bilinéaires symétriques semi-définies positives sur $(V/V'^\perp)\otimes \R$ dont le radical est défini sur~$\Q$. Si $V''\subset V'$, il existe une inclusion de $C(V/V''^\perp)$ dans $C(V/V'^\perp)$. Notons $\mathcal{C}$ le quotient de l'union disjointe
$$\coprod_{V'\in\mathfrak{C}} C(V/V'^\perp)\: ,$$
par la relation d'équivalence engendrée par les inclusions de $C(V/V''^\perp)$ dans $C(V/V'^\perp)$ pour tous sous-espaces totalement isotropes $V''\subset V'$. C'est un complexe conique qui est naturellement muni d'une action de $\gsp(\Z)$. Le sous-espace $\oplus_{j=g+1}^{2g} \Z\cdot x_j$ de $V$ est totalement isotrope et, si l'on fixe un isomorphisme $$X \isolong \bigoplus_{j=1}^{g} \Z\cdot x_j = V/ \left(\oplus_{j=g+1}^{2g} \Z\cdot x_j \right) ^\perp \: ,$$ on obtient une injection de $C(X)$ dans $\mathcal{C}$.
Cette injection est équivariante sous l'action de~$\mathrm{GL}(X)$, vu comme un sous-groupe de~$\gsp(\Z)$.

\subsubsection{Décompositions du complexe conique} On appelle \textit{décomposition polyédrale rationnelle de} $\mathcal{C}$ la donnée d'une décomposition polyédrale rationnelle de $C(V/V'^\perp)$ pour tout $V'\in\mathfrak{C}$~;~on demande en outre que la décomposition de $C(V/V''^\perp)$ soit la restriction à $C(V/V''^\perp)$ de la décomposition de $C(V/V'^\perp)$ si $V''\subset V'$. Soit $\mathfrak{S}$ une décomposition polyédrale rationnelle de $\mathcal{C}$.

\begin{definition3} On dit que $\mathfrak{S}$ est admissible relativement à~$\Gamma_0$ si elle est équivariante sous l'action de $\Gamma_0$ et si le quotient $\mathfrak{S}/\Gamma_0$ est fini. On dit que $\mathfrak{S}$ est lisse si pour tout  $V'\in\mathfrak{C}$, la décomposition induite par $\mathfrak{S}$ sur $C(V/V'^\perp)$ est lisse relativement à la structure entière duale de $S_{V',\:0}$.
\end{definition3}

\`A toute décomposition $\Sigma$ de $C(X)$ qui est $GL(X)$-admissible on peut associer son orbite sous~$\gsp(\Z)$ dans $\mathcal{C}$, ce qui détermine une décomposition~$\mathfrak{S}_\Sigma$ de $\mathcal{C}$ qui est $\Gamma_0$-admissible. \textit{A priori}, $\mathfrak{S}_\Sigma$ n'est pas lisse. En revanche, il existe une décomposition $\mathfrak{S}$ de $\mathcal{C}$ qui raffine $\mathfrak{S}_\Sigma$ et qui est lisse (\cite{Compact@AMRT}~II.5.3 ou  \cite{These@Pink}~9.20).

\subsubsection{Plongements toroïdaux associés à $\mathfrak{S}$} On fixe une décomposition $\Gamma_0$-admissible $\mathfrak{S}$ de $\mathcal{C}$. Elle pourra par exemple être égale à $\mathfrak{S}_\Sigma$, ou bien être lisse. Pour tout $\sigma \in \mathfrak{S}$, il existe un unique $V'_\sigma\in\mathfrak{C}$ tel que $\sigma$ soit inclus dans l'intérieur de $C(V/{V'_\sigma}^\perp)$. On dispose alors du champ algébrique $\mathcal{M}_{V'_\sigma,\: 0}$ et du plongement toroïdal 
$$  \mathcal{M}_{V'_\sigma,\: 0}  \hookrightarrow \overline{\mathcal{M}}_{V'_\sigma,\: \sigma,\: 0}$$
qui est construit à partir du cône $\sigma\subset \mathrm{Hom}(\mathrm{Sym}^2(V/{V'_\sigma}^\perp),\R)$ et de la structure entière $S_{V'_\sigma,\: 0}\subset \mathrm{Sym}^2(V/{V'_\sigma}^\perp)\otimes \R$.  Le champ algébrique $\overline{\mathcal{M}}_{V'_\sigma,\: \sigma,\: 0}$ est affine de type fini sur $\mathcal{B}_{V'_\sigma,\: 0}$. Il est muni d'une stratification localement fermée paramétrée par les faces de~$\sigma$. On note $Z_{V'_\sigma,\: \sigma,\: 0}$ l'unique strate fermée. 
Soit $V'\in\mathfrak{C}$. Notons
$$  \mathcal{M}_{V',\: 0}  \hookrightarrow \overline{\mathcal{M}}_{V',\: 0}$$
le plongement toroïdal localement de type fini sur $\mathcal{B}_{V',\: 0}$ associé à la restriction de $\mathfrak{S}$ à $C(V/V'^\perp)$, et à la structure entière $S_{V',\: 0}$. Le champ algébrique $\overline{\mathcal{M}}_{V',\: 0}$ est muni d'une stratification paramétrée par la restriction de $\mathfrak{S}$ à $C(V/V'^\perp)$. Notons $Z_{V',\: 0}$ l'union des strates paramétrées par des cônes inclus dans l'intérieur de $C(V/V'^\perp)$. Le champ $Z_{V',\: 0}$ est un sous-champ fermé de $\overline{\mathcal{M}}_{V',\: 0}$. 

\begin{remarque3} Si $\mathfrak{S}$ est lisse, les champs $\overline{\mathcal{M}}_{V'_\sigma,\: \sigma,\: 0}$ et $\overline{\mathcal{M}}_{V',\: 0}$ sont lisses sur $\mathcal{B}_{V',\: 0}$. Leurs singularités sont donc celles des champs~$\agrK$ pour $r\geq 0$. En particulier, les champs $\overline{\mathcal{M}}_{V'_\sigma,\: \sigma,\: 0}$ et $\overline{\mathcal{M}}_{V',\: 0}$ sont lisses sur $\Spec(\Z[1/p])$ et le complémentaire des strates ouvertes est un diviseur à croisements normaux.
\end{remarque3}

\subsubsection{Groupes d'automorphismes}  Soit $V'\in\mathfrak{C}$. Le stabilisateur de $V'$ dans $\Gamma_0$ agit sur $V/V'^\perp$. On note $\Gamma_{V'}$ son image dans $\mathrm{GL}(V/V'^\perp)$.
Le groupe $\Gamma_{V'}$ agit sur $C(V/V'^\perp)$ en stabilisant globalement la restriction de $\mathfrak{S}$ à $C(V/V'^\perp)$. Ce groupe agit donc sur $\overline{\mathcal{M}}_{V',\: 0}$. Le quotient $\overline{\mathcal{M}}_{V',\: 0}/\Gamma_{V'}$ est de type fini sur $\Spec(\Z)$ car $\mathfrak{S}/\Gamma_0$ est fini.
Pour tout cône $\sigma$ de~$\mathfrak{S}$, on note
$\Gamma_\sigma$ le stabilisateur de $\sigma$ dans ${\Gamma_{V'_\sigma}}$.
C'est un sous-groupe fini de $\mathrm{GL}(V/{V'_\sigma}^\perp)$ qui agit sur $\overline{\mathcal{M}}_{V'_\sigma,\: \sigma,\: 0}$.

\subsubsection{Cartes formelles associées à $\mathfrak{S}$}
Pour tout cône $\sigma$ de $\mathfrak{S}$, notons 
$$\widehat{\overline{\mathcal{M}}}_{V'_\sigma,\: \sigma,\: 0}$$
le complété formel de ${\overline{\mathcal{M}}}_{V'_\sigma,\: \sigma,\: 0}$ le long de $Z_{V'_\sigma,\: \sigma,\: 0}$. Le groupe fini $\Gamma_\sigma$ agit sur ce schéma formel, puisqu'il préserve la strate~$Z_{V'_\sigma,\sigma,\: 0}$.
Pour tout $V'\in\mathfrak{C}$, notons
$$\widehat{\overline{\mathcal{M}}}_{V',\: 0}$$
le complété formel de $\overline{\mathcal{M}}_{V',\: 0}$ le long de $Z_{V',\: 0}$. Il est muni d'une action de $\Gamma_{V'}$ et le quotient est formellement de type fini sur $\Spec(\Z)$.

\subsubsection{\'Enoncé du théorème principal}

Soient $\Sigma$ une décomposition de $C(X)$ lisse pour la structure entière $\mathrm{Sym}^2(X)$, et $\mathfrak{S}$ une décomposition de $\mathcal{C}$ qui raffine $\mathfrak{S}_\Sigma$. Notons $${\overline{\mathcal{A}}_g^{\Sigma}}$$ la compactification toroïdale de $\ag$ construite par Faltings et Chai.

\begin{theoreme3} \label{thPrinc} Il existe un champ algébrique $${\overline{\mathcal{A}}_{g,\: 0}^{\mathfrak{S}}}$$
qui est propre sur $\Spec(\Z)$, qui contient $\agK$ comme ouvert dense, et tel que l'on ait un diagramme cartésien 
\begin{eqnarray*}
\xymatrix{ \agK\ar[r] \ar[d] & {\overline{\mathcal{A}}_{g,\: 0}^{\mathfrak{S}}} \ar[d] \\
\ag \ar[r] & {\overline{\mathcal{A}}_g^{\Sigma}}. }
\end{eqnarray*}
Le champ ${\overline{\mathcal{A}}_{g,\: 0}^{\mathfrak{S}}}$ est muni d'une stratification localement fermée paramétrée par $\mathfrak{S}/\Gamma_0$. 
Son complété formel le long de la strate associée à tout cône $\sigma$ de $\mathfrak{S}$ est isomorphe à
$$\widehat{\overline{\mathcal{M}}}_{V'_\sigma,\: \sigma,\: 0} \: /\: \Gamma_\sigma\: .$$
Pour tout $V'\in\mathfrak{C}$, le complété formel de ${\overline{\mathcal{A}}_{g,\: 0}^{\mathfrak{S}}}$ le long de l'union des strates associées aux cônes de $\mathfrak{S}$ qui sont inclus dans l'intérieur de $C(V/V'^\perp)$ est isomorphe à
$$\widehat{\overline{\mathcal{M}}}_{V',\: 0} \: /\: \Gamma_{V'}\: .$$
Le champ algébrique ${\overline{\mathcal{A}}_{g,\: 0}^{\mathfrak{S}}}$ est normal et sa projection  sur $\Spec(\Z)$ est plate à fibres géométriques réduites.
Si $\mathfrak{S}$ est lisse, ses singularités sont celles des champs $\agrK$ pour $r\geq 0$. En particulier, il est lisse sur $\Spec(\Z[1/p])$ et le complémentaire de $\agK$ est un diviseur relatif à croisements normaux.
\end{theoreme3}

Dans le cas où $\mathfrak{S}=\mathfrak{S}_\Sigma$, on a ${\overline{\mathcal{A}}_{g,\: 0}^{\Sigma}}={\overline{\mathcal{A}}_{g,\: 0}^{\mathfrak{S}}}$ et le théorème~\ref{thPrinc} précise la combinatoire évoquée après le théorème~\ref{thPrinc1}.

\medskip

\begin{demo} On reprend la construction de ${\overline{\mathcal{A}}_{g,\: 0}^{\Sigma}}$ en l'adaptant légèrement. Pour tout $\sigma\in\Sigma$, notons $\mathfrak{S}_\sigma$ l'ensemble des cônes de $\mathfrak{S}$ qui appartiennent à l'orbite de $\sigma$ sous l'action de $\gsp(\Z)$. Posons 
$$\overline{\mathcal{M}}_{\mathfrak{S}_\sigma,\:0} \: = \: \coprod_{\tau\in \mathfrak{S}_\sigma} \overline{\mathcal{M}}_{V'_\tau,\tau,\: 0}\: .$$
Soit $\overline{\mathcal{M}}_\sigma^{et}$ une carte locale algébrique sans niveau. On pose
$$\overline{\mathcal{M}}_{\mathfrak{S}_\sigma,\:0}^{et} \: = \: \overline{\mathcal{M}}_{\mathfrak{S}_\sigma,\:0} \times_{\overline{\mathcal{M}}_{_\sigma}} \overline{\mathcal{M}}_{\sigma}^{et}\: ,$$
On définit alors $\overline{U}_{\mathfrak{S}, \:0}$ comme l'union de tous les schémas $\overline{\mathcal{M}}_{\mathfrak{S}_\sigma,\:0}^{et}$ pour $\sigma \in \Sigma$. On définit le groupoïde $\overline{R}_{\mathfrak{S},\: 0}$ comme dans le paragraphe~\ref{partieRelEqNiveau} et l'on vérifie qu'il est étale. On note
$${\overline{\mathcal{A}}_{g,\: 0}^{\mathfrak{S}}} \: = \: \overline{U}_{\mathfrak{S}, \:0} \: / \: \overline{R}_{\mathfrak{S},\: 0} $$
le champ algébrique quotient. On voit dans les démonstrations des propositions analogues à~\ref{propRbEt} et~\ref{propR0bet} que le groupoïde $\overline{R}_{\mathfrak{S},\: 0}$ respecte la stratification et qu'il identifie des  strates conjuguées sous l'action de $\Gamma_0$. On en déduit les assertions relatives à la stratification de la compactification. Les autres assertions résultent de~\cite{FlatSympl@Gortz} et de la proposition~\ref{propAgKNorm}.
\end{demo}

\subsubsection{Propriété modulaire à la frontière}

Soient $(U\hookrightarrow S)$ un objet de $\mathsf{Ouv}$ et $(G',H'_\bullet)$ un objet de $\agK(U)$ tel que $G'$ s'étend en un schéma semi-abélien sur $S$. Notons $G'_i = G'/H'_i$ pour $1\leq i \leq s$.  Il existe une chaîne d'isogénies 
$$G'=G'_0 \: \longrightarrow  \: G'_1 \: \longrightarrow \: G'_2 \: \longrightarrow \: \cdots \: \longrightarrow \: G'_s $$
entre schémas semi-abéliens sur $S$. Rappelons qu'on a défini le dual d'un schéma semi-abélien dans le paragraphe~\ref{parGlobB}. Il existe donc une chaîne de morphismes
\begin{equation}\label{eqchaine}
\underline{X}(G'_0\: ^t) \longrightarrow \underline{X}\left(G'_1\: ^t\right) \longrightarrow \cdots \longrightarrow \underline{X}\left(G'_s\: ^t\right) \longrightarrow \underline{X}\left(G'_s \: \right) \longrightarrow \cdots \rightarrow \underline{X}(G'_0)\: .
\end{equation}
Supposons qu'il existe $V'\in\mathfrak{C}$ et une surjection de $V/V'^\perp$ dans $\underline{X}(G)$ qui identifie l'image de $\mathbb{V}_\mathcal{D}^\bullet$ et~(\ref{eqchaine}). D'après la proposition~\ref{propXi0w}, la famille des formes bilinéaires $\underline{B}(G'_i)$ définit un morphisme $\underline{B}(G',H'_\bullet)$ de $S_{V',\: 0}$ dans $\underline{\mathrm{Div}}_S$. La propriété modulaire à la frontière suivante est la traduction de la proposition~\ref{propModFrontAgKbVersion1} en terme de complexe conique.

\begin{proposition2} \label{propModFrontAgKb} Le groupoïde $\left(\agK\hookrightarrow {\overline{\mathcal{A}}_{g,\: 0}^{\mathfrak{S}}}\right)(U\hookrightarrow S)$
paramètre les objets $(G',H'_\bullet)$ de $\agK(U)$ tels que $G'$ s'étend en un schéma semi-abélien sur $S$, et tels que localement pour la topologie étale de $S$, il existe $V'\in\mathfrak{C}$, une surjection $V/V'^\perp\rightarrow \underline{X}(G')$ identifiant l'image de $\mathbb{V}_\mathcal{D}^\bullet$ et~(\ref{eqchaine}), et  $\sigma \in \mathfrak{S}$ dont l'adhérence contient la forme bilinéaire $p_D \circ \underline{B}(G',H'_\bullet)$ pour tout diviseur irréductible réduit $D$ de $S$.
\end{proposition2}

\section{Application : existence du sous-groupe canonique}

Le champ $\agK$ associé à $\mathcal{D}=\{g\}$ paramètre les variétés abéliennes principalement polarisées munies d'un sous-groupe lagrangien de $p$-torsion. Notons $\agK^{\mathrm{ord-mult}}$ le lieu ordinaire-multiplicatif de $\agK\times\Spec(\F_p)$~;~c'est par définition le sous-schéma ouvert de $\agK\times \Spec(\F_p)$ sur lequel le schéma abélien universel est ordinaire et où le sous-groupe lagrangien universel est égal au noyau du morphisme de Frobenius. Notons de même $\ag^{\mathrm{ord}}$ le lieu ordinaire de $\ag\times \Spec(\F_p)$. L'oubli du niveau induit un isomorphisme $$\agK^{\mathrm{ord-mult}}\:\isolong\:\ag^{\mathrm{ord}}$$
dont l'application réciproque associe à une variété abélienne ordinaire le noyau du morphisme de Frobenius. Nous étendons ce résultat aux compactifications toroïdales, et en déduisons une nouvelle démonstration de l'existence du sous-groupe canonique pour des familles de variétés abéliennes (\cite{Canonique@AbbesMokrane} et~\cite{Canonique@Andreatta}), selon une idée de Vincent Pilloni.

\subsection{Lieu ordinaire des compactifications} \label{partieOrd}

Soit $\Sigma$ une décomposition polyédrale admissible de $C(X)$. On note $\agKb$ et $\agb$ les  compactifications toroïdales de $\agK$ et $\ag$ qui lui sont associées. D'après le théorème~\ref{thPrinc1}, il existe un morphisme de $\agKb$ dans $\agb$ qui étend l'oubli du niveau. Soit $S$ le spectre d'un corps de caractéristique $p$.

\begin{definition2}  Un schéma semi-abélien sur $S$ est ordinaire s'il est extension d'une variété abélienne ordinaire par un tore.
\end{definition2}

Notons $G$ le schéma semi-abélien canonique sur $\agb$ et $H$ le sous-groupe quasi-fini et plat canonique sur $\agKb$. Notons $$\agKb^{\mathrm{ord-mult}}$$ le sous-schéma ouvert de $\agKb\times \Spec(\F_p)$ sur lequel $G$ est ordinaire et $H$ est égal au noyau du morphisme de Frobenius $\mathrm{Frob}_G$ de $G$. Notons de même $$\agb^{\mathrm{ord}}$$ le sous-schéma ouvert de $\agb\times \Spec(\F_p)$ sur lequel $G$ est ordinaire.

\begin{remarque2} Soit $\sigma\in\mathfrak{S}_\Sigma$. Réutilisons les notations des paragraphes~\ref{partieCompactAmeliore} et~\ref{partieCompactAmeliore}. Notons $f$ l'isogénie universelle entre $1$-motifs universels sur $\mathcal{M}_{{V'_\sigma},\: \sigma,\: 0}$. Le sous-schéma $$\agKb^{\mathrm{ord-mult}}$$ rencontre la $\sigma$-strate de $\agKb$ si et seulement si le noyau de $f$ est égal à partie semi-abélienne, ce qui signifie $W_1(\mathrm{Ker}(f))=\mathrm{Ker}(f)$. On peut reformuler cette condition en disant que $\mathbb{V}_\mathcal{D}^0$ et $\mathbb{V}_\mathcal{D}^1$ ont même image dans $V/{V'_\sigma}^\perp$.
\end{remarque2}

La catégorie fibrée à la frontière $$(\agK^{\mathrm{ord-mult}}\hookrightarrow  \agKb^{\mathrm{ord-mult}})$$ vérifie une propriété modulaire analogue à celle de $(\agK\hookrightarrow \agKb)$, donnée par la proposition~\ref{propModFrontAgKb}. Soit $(U\hookrightarrow S)$ un objet de $\mathsf{Ouv}$ tel que $S$ soit un schéma sur $\Spec(\F_p)$.

\begin{proposition2} Le groupoïde $\left(\agK^{\mathrm{ord-mult}}\hookrightarrow  \agKb^{\mathrm{ord-mult}} \right)(U\hookrightarrow S)$
paramètre les objets $G'$ de $\ag^{\mathrm{ord}}(U)$ tels que $G'$ s'étend en un schéma semi-abélien ordinaire sur $S$, et tels que localement pour la topologie étale de $S$, il existe
\begin{itemize}
\item $V'\in\mathfrak{C}$ tel que $\mathbb{V}_\mathcal{D}^0/V'^\perp=\mathbb{V}_\mathcal{D}^1/V'^\perp$, 
\item une surjection $V/V'^\perp\rightarrow \underline{X}(G')$,
\item $\sigma \in \mathfrak{S}_\Sigma$ dont l'adhérence contient la forme bilinéaire $p_D \circ \underline{B}(G',\mathrm{Ker}(\mathrm{Frob}_{G'}))$ pour tout diviseur irréductible réduit $D$ de $S$.
\end{itemize}
\end{proposition2}

Par conséquent, le morphisme de $\ag^{\mathrm{ord}}$ dans $\agK^{\mathrm{ord-mult}}$ induit par le noyau du morphisme de Frobenius s'étend en un morphisme $$\agb^{\mathrm{ord}}\rightarrow \agKb^{\mathrm{ord-mult}}\: .$$
Par densité de $\agK^{\mathrm{ord-mult}}$ dans $\agKb^{\mathrm{ord-mult}}$ et de $\ag^{\mathrm{ord}}$ dans $\agb^{\mathrm{ord}}$, ce morphisme est l'inverse du morphisme d'oubli du niveau. Ainsi, les schémas $\agKb^{\mathrm{ord-mult}}$ et $\agb^{\mathrm{ord}}$ sont isomorphes.

\subsection{Surconvergence du tube}

Notons $\agb^{\mathrm{for}}$ le complété formel de $\agb$ le long de $\agb\times \Spec(\F_p)$, notons  $\agb^{\mathrm{rig}}$ la fibre générique rigide du schéma formel $p$-adique $\agb^{\mathrm{for}}$, et notons $T$ le tube de $\agb^{\mathrm{ord}}$ dans $\agb^{\mathrm{rig}}$. Enfin, notons $G$ le schéma semi-abélien canonique sur $\agb^{\mathrm{rig}}$.

\begin{theoreme2} \label{theo} Il existe un voisinage strict $U$ de $T$ dans $\agb^{\mathrm{rig}}$ et un sous-groupe fini, plat et lagrangien $H$ de $G_U[p]$ tel que $H_T$ soit le relèvement du noyau du morphisme de Frobenius.
\end{theoreme2}

\begin{demo}
Les arguments de la partie~\ref{partieOrd} montrent l'existence d'un diagramme commutatif
$$\xymatrix{
& \agKb \times \Spec(\F_p) \ar[d] \ar[r] & \agKb^{\mathrm{for}} \ar[d]^u \\ 
\agb^{\mathrm{ord}} \ar[ru] \ar[r] & \agb \times \Spec(\F_p) \ar[r] &\agb^{\mathrm{for}} }$$
où $\agKb^{\mathrm{for}}$ désigne le complété formel de $\agKb$ le long de $\agKb\times\Spec(\F_p)$. Le morphisme $u$ est propre. D'après~\cite[lem.~4.18]{DegCourbe@Mumford} et~\cite{FlatSympl@Gortz}, il réalise un isomorphisme sur un voisinage de $\agb^{\mathrm{ord}}$. En particulier, $u$ est étale sur un voisinage de $\agb^{\mathrm{ord}}$.
On note $\agKb^{\mathrm{rig}}$ l'espace rigide associé à $\agKb^{\mathrm{for}}$ et $T_0$  le tube de $\agb^{\mathrm{ord}}$ dans $\agKb^{\mathrm{rig}}$. En appliquant~\cite[th.~1.3.5]{Rigide@Berthelot}, on obtient l'existence d'un voisinage strict $U_0$ (resp. $U$) de $T_0$ (resp. de $T$) dans $\agKb^{\mathrm{rig}}$ (resp. $\agb^{\mathrm{rig}}$) tels que $u$ induise un isomorphisme
$$u \: : \: U_0 \: \isolong U \: .$$
L'application réciproque définit le sous-groupe lagrangien recherché sur $U$.
\end{demo}

Ainsi, les variétés abéliennes paramétrées par $U\cap \ag^{\mathrm{rig}}$ admettent un sous-groupe canonique. On obtient un résultat moins précis que ceux de~\cite{Canonique@AbbesMokrane} et~\cite{Canonique@Andreatta}, car il ne donne pas de condition suffisante pour qu'une variété abélienne admette un sous-groupe canonique. Toutefois, il suffit pour la plupart des applications aux formes modulaires de Siegel $p$-adiques.

\newpage

\bibliography{bibliographie}
\bibliographystyle{smfalpha}

\end{document}

%% file: paquets.tex
\usepackage{amsmath,a4wide}
\usepackage{amssymb}
\usepackage{xspace}
\usepackage{euscript}
\usepackage{mathrsfs} 
\usepackage[all]{xy}
\usepackage{hyperref}
\usepackage[francais]{babel}
\usepackage[T1]{fontenc}
\usepackage{textcomp}
\usepackage{smfthm}

%% file: numerotation.tex
\renewcommand{\theenumi}{\textbf{\roman{enumi}}} 

\renewcommand\theequation{\thesection.\textbf{\alph{equation}}} 

\makeatletter
\@addtoreset{equation}{subsection}
\makeatother

\renewcommand\thesubsubsection{\textbf{\thesubsection.\arabic{subsubsection}}}
\renewcommand\thesubsection{\textbf{\thesection.\arabic{subsection}}}

%% file: lemme-proposition-theoreme.tex
\newtheorem{lemme2}[subsubsection]{Lemme}	
\newtheorem{lemme3}[paragraph]{Lemme}		

\newtheorem{proposition2}[subsubsection]{Proposition}
\newtheorem{proposition3}[paragraph]{Proposition}

\newtheorem*{thprinc}{Théorème principal}
\newtheorem*{thrien}{Théorème}

\newtheorem{theoreme2}[subsubsection]{Th\'eor\`eme}
\newtheorem{theoreme3}[paragraph]{Th\'eor\`eme}

\newtheorem{corollaire3}[paragraph]{Corollaire}

\theoremstyle{definition}

\newtheorem{definition2}[subsubsection]{D\'efinition}
\newtheorem{definition3}[paragraph]{D\'efinition}

\newtheorem{remarque2}[subsubsection]{Remarque}
\newtheorem{remarque3}[paragraph]{Remarque}

\newenvironment{demo}{\noindent{\textit{D\'emonstration. --- }}}{~\qedsymbol \vspace{4mm}}